\pdfoutput=1

\documentclass[preprint, 12pt]{elsarticle}
\usepackage[left=2cm,right=1.5cm, top=2.5cm,bottom=2.5cm]{geometry}

\usepackage{amssymb}

\usepackage{bm}
\usepackage{amsmath}
\usepackage{mathtools}
\usepackage{xcolor}
\usepackage{algorithm,algpseudocode}
\usepackage{siunitx}
\usepackage{arydshln}
\usepackage{empheq}
\usepackage{subcaption}
\usepackage{adjustbox}
\usepackage{stmaryrd}
\usepackage{amsfonts}
\usepackage{cases}
\usepackage{hyperref}

\newcommand{\norm}[1]{\left\lVert#1\right\rVert}
\newcommand{\fe}[1]{\mathbf{#1}}
\newcommand{\pr}[1]{\left(#1\right)}
\newcommand{\R}{\mathbb{R}}
\newcommand{\Ste}{^{\rm s}}
\newcommand{\Con}{^{\rm c}}

\newcommand{\stressPSte}{{\rm N}_\mu}
\newcommand{\stressPCon}{\sigma_\mu}
\newcommand{\strainP}{\varepsilon_\mu}
\newcommand{\stressP}{\sigma_\mu}
\newcommand{\intvarP}{\gamma_\mu}
\newcommand{\IntVarEqua}{\mathcal{F}^{\gamma}_\mu}
\newcommand{\IntVarEquaCon}{\mathcal{F}^{\gamma^{\rm c}}_\mu}
\newcommand{\IntVarEquaSte}{\mathcal{F}^{\gamma^{\rm s}}_{\mu}}

\newcommand{\FEstressP}{\bm{\sigma}_\mu}

\newcommand{\FEdisP}{\mathbf{u}_\mu}

\newcommand{\HFresDis}{\mathbf{R}^{\rm hf}_\mu}

\newcommand{\ConsEquaSte}{\mathcal{F}^{\rm N}_{\mu}}
\newcommand{\StepQSOp}{\mathcal{F}_\mu}

\journal{Elsevier}

\begin{document}

\begin{frontmatter}


\title{Projection-based model order reduction for prestressed concrete with an application to the standard section of a nuclear containment building}

\author[1,2,3]{Eki Agouzal\corref{cor1}}\ead{eki.agouzal@inria.fr}
\author[1]{Jean-Philippe Argaud}
\author[2,3]{Michel Bergmann}
\author[1]{Guilhem Ferté}
\author[1]{Sylvie Michel-Ponnelle}
\author[2,3]{Tommaso Taddei}
\cortext[cor1]{Corresponding author}

\affiliation[1]{organization={EDF Lab Paris-Saclay},
            addressline={7 Boulevard Gaspard Monge}, 
            city={Palaiseau},
            postcode={91120}, 
            state={},
            country={France}}

\affiliation[2]{organization={IMB, UMR 5251, Univ. Bordeaux},
            addressline={}, 
            city={Talence},
            postcode={33400}, 
            state={},
            country={France}}
            
\affiliation[3]{organization={INRIA, Inria Bordeaux Sud-Ouest, Team MEMPHIS, Univ. Bordeaux},
            addressline={}, 
            city={Talence},
            postcode={33400}, 
            state={},
            country={France}}           

\begin{abstract}\footnotesize
We propose a projection-based model order reduction procedure for the ageing of large prestressed concrete structures. Our work is motivated by applications in the nuclear industry, particularly in the simulation of containment buildings. Such numerical simulations involve a multi-modeling approach: a three-dimensional nonlinear thermo-hydro-visco-elastic rheological model is used for concrete; and prestressing cables are described by a one-dimensional linear thermo-elastic behavior. A kinematic linkage is performed in order to connect the concrete nodes and the steel nodes: coincident points in each material are assumed to have the same displacement. We develop an adaptive algorithm based on a Proper Orthogonal Decomposition (POD) in time and greedy in parameter to build a reduced order model (ROM). The nonlinearity of the operator entails that the computational cost of the ROM assembly scales with the size of the high-fidelity model. We develop an hyper-reduction strategy based on empirical quadrature to bypass this computational bottleneck: our approach relies on the construction of a reduced mesh to speed up online assembly costs of the ROM. We provide numerical results for a standard section of a double-walled containment building using a qualified and broadly-used industrial grade finite element solver for structural mechanics (code$\_$aster).
\end{abstract}


\end{frontmatter}


\section{Introduction}
\subsection{Context}

A nuclear power plant is an industrial facility designed to produce electricity, and whose nuclear steam supply comprises one or more nuclear reactors. Électricité De France (EDF) operates a fleet of 56 reactors, 24 of which have so-called double-walled nuclear containments buildings (NCBs). In this case, the safety of the nuclear plants rely on an outer wall made of reinforced and prestressed concrete that shield the reactor form external aggression and a inner wall made of prestressed concrete (no steel liner) that should contain any leaks of radioelements in case of accident. However, the leakage rate may be influenced by the ageing of these large concrete structures. This phenomenon is mainly due to two physical phenomena: drying and creep of concrete.  Creep and drying induce delayed strains, and, thus, a loss of prestressing effects. All these phenomena may lead to a modification of the concrete's permeability, or to the re-opening of cracks within the material. These changes can result in an increase in the leakage rate through the concrete. Therefore, the mechanical response of the inner wall is carefully monitored using a set of deformation sensors embedded in the concrete, and the leak-tightness of the inner containment is checked every 10 years thanks to A Integrated Leakage Rate Test, during which the NCB's internal relative pressure rises to 4.2 bars. These inspections play a crucial role in ensuring that the structure maintains its optimal operational condition.\\

In recent years, research has been carried out into the realistic modeling of the thermo-hydro-mechanical (THM) behavior, and even leakage (THM-L), of concrete in large prestressed concrete structures. In view of the complexity of the phenomena involved in modeling these structures, various techniques may be applied. More specifically, existing numerical approaches in the literature for modeling concrete ageing can be divided into two main categories: strong coupling strategies \cite{dal2007multiphase}\cite{gawin2003modelling}, where all dependencies between behaviors are accounted for, or weak coupling strategies \cite{asali2016numerical}\cite{bouhjiti2018accounting}\cite{jason2007hydraulic} (chained calculations) which aim to reduce these inter-dependencies by neglecting, e.g., the effect of mechanical stresses on thermal and hydric responses. The aim of these numerical models is to predict the temporal behavior of physical quantities of interest (QoIs), such as water saturation in concrete, delayed deformations, and stresses. Comprehensive understanding of these diverse fields has facilitated the development of numerical methods for estimating leakage rates, notably utilizing prestress loss in cables \cite{charpin2022predicting}. \\

Achieving accurate simulations for NCB systems involves handling a potentially large number of model parameters, often with limited available knowledge. As noted in \cite{bouhjiti2020stochastic}, numerous parameters lack sufficient information, leading to the need for expert judgment in quantifying uncertainties \cite{de2008uncertainty}. The uncertainties in the output fields of numerical calculations are hence significant and might be linked to the inadequacy of the PDE model (structural uncertainty) or to the calibration of the parameters (parametric uncertainty). To address this issue, auscultation data, which are obtained for studying the long-term behavior of the structure, offer valuable insights. Those data, provided by the monitoring structures, can be leveraged to further enhance understanding of the system’s response. \\

The past decade has witnessed significant progress in the development of numerical methods that combine data and models --- in effect data assimilation --- for THM systems. Bayesian inference has been applied as a first step to predict the THM-L behavior of confinement structures. To reduce the computational burden, Bayesian methods have been implemented in combination with simplified models of the system response: \citet{berveiller2012updating} employed Bayesian inference to refine predictions of deformations using a simplified one-dimensional model of NCB; in a more recent study \cite{rossat2021bayesian}, Bayesian updating of NCB leak response was presented, based on a simplified one-dimensional model. On the other hand, \citet{rossat2022bayesian}  extended Bayesian strategies to three-dimensional models, employing a 1:3 NCB \cite{masson2014vercors} with a metamodel founded on a finite element model of a representative structural volume (RSV). In addition to Bayesian approaches, other methodologies are deployed to address uncertainties in parameters. Variational assimilation methodologies (3D VAR) are utilized to integrate a priori knowledge of parameters with observations, providing an alternative strategy to address model uncertainties.

\begin{figure}[h!]
\begin{center}
\includegraphics[scale=0.95]{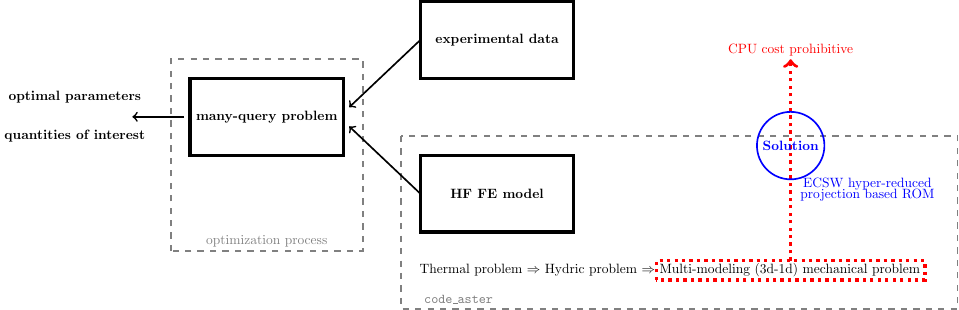} 
\end{center}
\caption{Highlighting the relevance of the ROM methodology for an industrial application: application to an HF finite element (FE) model for the simulation of a containment building modeled by a THM approach for prestressed concrete \label{sec:intro:fig:context}}
\end{figure}

The data assimilation strategies discussed so far pose challenges as they involve solving many-query problems. For practical high-fidelity (HF) models, data assimilation strategies result in prohibitive costs, which partly explains the scarcity of reported results on three-dimensional models. FE simulations are key to achieve detailed estimations of the temporal behavior of THM QoIs; however, the simulation of the aging of a real containment building over several decades takes approximately a whole day, even with parallel computations, and is hence impractical for many-query scenarios. \\

Parametric Model Reduction (pMOR) is a family of algorithms aimed at reducing the marginal cost associated with the computation of the solution to a parametric problem. This reduction is achieved by leveraging prior knowledge obtained from previously conducted HF calculations, allowing for the approximation of a field over a range of parameters. Our objective is to develop an intrusive pMOR procedure for the mechanical simulation of double-walled power plant containment buildings; we consider the application to a NCB, which involves a complex FE model of a RSV (cf. Figure \ref{sec:intro:fig:context}).

\subsection{Objective of the paper and relation to previous works}

The key contribution of this paper is the development of a hyper-reduced model for nonlinear mechanics problems within the multi-modeling framework. More specifically, we develop an approach that provides a high-quality reduced order model (ROM \cite{hesthaven2016certified}\cite{quarteroni2015reduced}\cite{rozza2008reduced}) to mimic the behavior of prestressed concrete, with an application to a standard section of a NCB . Clearly, such a problem falls within the scope outlined above: it consists of concrete in which tendons are embedded. Each material has its own constitutive equation, and the two are kinematically coupled. Furthermore, our approach aims to develop a ROM useful for engineering applications, which means that in addition to being able to approximate the solution, it must provide QoIs close to those obtained for the HF solution: prestress loss in the cables and tangential and vertical deformations inside and outside the standard section.  \\

Our strategy for building a ROM is founded on previous work. The methodology relies on a Galerkin projection method. We develop an adaptive algorithm based on a Proper Orthogonal Decomposition (POD)-Greedy strategy \cite{haasdonk2008reduced} to construct a ROM. This algorithm is an approach that iteratively improves the model using poorly-approximated solutions, so as to get a reduced model valid on a set of parameters. We rely on the Proper Orthogonal Decomposition (POD \cite{berkooz1993proper}, \cite{bergmann2009enablers}, \cite{volkwein2011model}) to compress the temporal trajectory of the physical problem. The nonlinearity of the operator entails that the computational cost of the ROM assembly scales with the size of the HF model. We develop an hyper-reduction \cite{hesthaven2016empirical}\cite{carlberg2013gnat}\cite{ryckelynck2005priori}\cite{hernandez2017dimensional} strategy based on empirical quadrature (EQ \cite{farhat2014dimensional}, \cite{farhat2015structure}) to bypass this computational bottleneck: our approach relies on the construction of a reduced mesh to speed up online assembly costs of the ROM. \\

The methodology developed and the simulations were carried out with a qualified and broadly-used industrial grade finite element solver for structural mechanics (\textsf{code$\_$aster}) \cite{aster}. Our work constitutes a continuation of research efforts at EDF R\&D to develop approaches for nonlinear mechanical problems in structural mechanics, with the aim of simulating real-world problems. In this respect, we mention previous work on nonlinear parabolic thermo-mechanical problems \cite{benaceur2018reduction}, on vibro-acoustic problems \cite{khoun2021reduced} and also on welding \cite{dinh2018modeles}. The approaches developed here must take into account the \textsf{code$\_$aster} computational framework, namely the dualization of the boundary conditions. In particular, the work in Reference \cite{dinh2018modeles} features one of the first efforts to design hyper-reduced ROMs in \textsf{code$\_$aster}. Our point of departure is the pMOR methodology of [1] that relies on Energy-Conserving Sampling and Weighting method (ECSW \cite{farhat2014dimensional}) for hyper-reduction to deal with a three-dimensional elasto-plastic holed plate.

\subsection{Layout of the paper}

The outline of the paper is as follows. In section \ref{sec:ROMmethodo}, we present the pMOR methodology. In section \ref{sec:THMmodel}, we present the formulation of the multi-modeling framework for the THM modelisation of prestressed concrete. In section \ref{sec:numres_srpb}, we validate the methodology for a non-parametric problem, before presenting numerical results for a parametric problem in section \ref{sec:numres_param}.

\section{Methodology for the ROM for the multi-modeling nonlinear mechanical problem\label{sec:ROMmethodo}}
\subsection{Formulation of the nonlinear quasi-static multi-modeling mechanical problem\label{sec:ROMmethodo:subsec:formqs}}
\subsubsection{Continuous formulation of the problem\label{sec:ROMmethodo:subsec:formqs;subsubsec:continuous}}

In this contribution, we study quasi-static nonlinear problems for mechanics. We focus on small-strain small displacements problems. We consider the modeling of large prestressed concrete structures. Therefore, the developed mechanical model is built on a coupling between a three-dimensional model (modeling the concrete) and a one-dimensional model (modeling the prestressing steel cables). We consider a domain $\Omega\subset \R^3$ of the space supposed to be sufficiently regular. As mentioned above, we assume that the domain $\Omega$ can be split into a three-dimensional domain $\Omega^{\rm c}$, and a one-dimensional domain $\Omega^{\rm s}$. The latter can be decomposed in $n_\mathcal{C}$ cables $\Omega^{\rm s}=\left\{\mathcal{C}_i\right\}_{i=1}^{n_\mathcal{C}}$, modeled by curves that correspond to their mean line. We introduce a vector of parameters $\mu\in \mathcal{P}\subset \R^p$, which contains physical parameters of the problem (coefficients in the constitutive equations of the steel or the concrete).\\

We denote by $u_\mu$ the vector of displacements, whether in cables or concrete and we denote by $\mathcal{X}$ the Hilbert space to which the field $u_\mu$ belongs. To identify the displacements in each subdomain, we shall note $u_\mu\Con$ the displacement in the concrete, and $u_\mu\Ste$ the displacement in the steel. Both of those fields can be seen as restrictions of $u_\mu$ on the corresponding domain. The mechanical strains tensor within the concrete is the symmetric gradient of the displacement and is denoted $\varepsilon^{c}_\mu = \nabla_s u_\mu\Con = \frac{1}{2}\left( \nabla u_\mu^{\rm c} + (\nabla u_\mu^{\rm c}\right)^\top)$, and the strains within the cables (also called uniaxial strains) are defined as $\strainP\Ste = \partial_s u_\mu\Ste$, where $\partial_s(.)$ is the derivative along the cable. We denote the stress tensor within the concrete $\stressPCon$, the normal forces in the steel $\stressPSte$, and the internal variables in the concrete $\gamma_\mu\Con$ and in the steel $\gamma_\mu\Ste$. We assume that the constitutive equations used depend on auxiliary variables, which we shall refer to in this section as a vector $H$. The fields enclosed in $H$ include previously computed fields and solutions to PDEs that do not depend on the parameters set in the vector $\mu$. The vector is comprised of fields that may appear and be used in the problem's constitutive or evolution equations. In the application case presented, namely in the case of a thermo-hydro-activated mechanical problem, this vector consists of the pair made of temperature and water content in the concrete. Details are provided in section \ref{sec:THMmodel}. We introduce the quasi-static equilibrium equations for the three-dimensional model, where we omit to specify the initial conditions (ICs) and the boundary conditions (BCs) for each subdomains:

$$\left\{
\begin{array}{rcl}
- \nabla \cdot \stressP &=& f_{\rm c} \quad \text{on} \ \Omega\Con,\\ \\
\stressPCon &=& \mathcal{F}_\mu^\sigma \left(\strainP\Con, \ \intvarP\Con, \ H\right), \\ \\
\dot{\gamma}_\mu\Con &=& \IntVarEquaCon\left(\stressPCon, \intvarP\Con, \ H\right),
\end{array}
\right.\quad \text{and} \quad
\left\{
\begin{array}{rcl}
\displaystyle \frac{\partial \stressPSte}{\partial s}&=& f_{\rm s}\quad \text{on} \ \Omega\Ste,\\ \\
\stressPSte &=& \ConsEquaSte\left(\partial_s u_\mu\Ste, \ \intvarP\Ste, \ H\right), \\ \\
\dot{\gamma}_\mu\Ste &=& \IntVarEquaSte\left(\text{N}_\mu, \intvarP\Ste, \ H\right),
\end{array}
\right.$$

\noindent where $\mathcal{F}_\mu^\sigma$ (resp. $\ConsEquaSte$) stands for the constitutive equation for the three-dimensional (resp. one-dimensional) problem, while the nonlinear operator $\IntVarEquaCon$ (resp. $\IntVarEquaSte$) denotes an equation of evolution of internal variables within the concrete (resp. the steel). To provide more compact notations, we introduce new notations for the fields defined on the whole domain, namely for the displacements, strains, generalized forces (stresses or normal efforts), internal variables and the loadings. All the details are provided in Table \ref{sec:ROMmethodo:subsec:formulation:table:notations}.

\begin{table}[h!]
\begin{center}
\begin{tabular}{cccc}
\hline
\textbf{Notation on $\Omega$} & \textbf{Notation on $\Omega^{\rm s}$} & \textbf{Notation on $\Omega^{\rm c}$} & \textbf{Definition} \\ \hline
$\mathfrak{S}_\mu$ & $\stressPSte$ & $\stressPCon$ & Generalized force \\ 
$u_\mu$ & $u_\mu^{\rm s}$ & $u_\mu^{\rm c}$ & Displacement \\  
$\varepsilon_\mu$ & $\varepsilon_\mu^{\rm s} = \partial_s u_\mu\Ste$ & $\varepsilon_\mu^{\rm c}$ & Strain \\ 
$\gamma_\mu$ & $\gamma_\mu^{\rm s}$ & $\gamma_\mu^{\rm c}$ & Internal variables \\
$f$ & $f_{\rm c}$ & $f_{\rm s}$ & External loading \\
\end{tabular}
\end{center}
\caption{Notations of the fields defined on the whole computational domain $\Omega$, whose definition depends on the subdomains ($\Omega^{\rm c}$ or $\Omega^{\rm s}$) \label{sec:ROMmethodo:subsec:formulation:table:notations}}
\end{table}

These notations enable us to recast the problem in a compact form, which helps to manage the multi-modeling (3d-1d) using three operators, $\mathcal{G}_\mu\left(.\right)$ for the equilibrium equation, $\mathcal{F}_\mu^\mathfrak{S}\left(.\right)$ for the constitutive equation and $\IntVarEqua\left(.\right)$ for the evolution equation for internal variables:

\begin{equation*}
\left\{
\begin{array}{rcl}
\mathcal{G}_\mu\left(\mathfrak{S}_\mu\right) &=& f,\\
\mathfrak{S}_\mu &=& \mathcal{F}_\mu^{\mathfrak{S}}\left(\mathfrak{S}_\mu, \ \intvarP, \ H\right), \\
\dot{\gamma}_\mu &=& \IntVarEqua\left(\varepsilon_\mu, \ \intvarP, \ H\right),
\end{array}
\right.
\end{equation*}

\noindent where we still omit the ICs and BCs used. In our study, the initial state of the problem is the material at rest, so all physical fields are assumed to be zero initially. The temporal discretization of the equations is done using a one-step integrator ($u_\mu^{(k+1)} = u_\mu^{(k)} + \Delta u_\mu^{(k+1)}$), which implies that the knowledge of the mechanical state is derived from the mechanical state previously computed (and the knowledge of the field $H$ at the current time). In our study, we consider both non-homogeneous Neumann conditions (defined on $\Gamma_{\rm n}^{\rm c}$ for the concrete) and homogeneous Dirichlet conditions for suitable linear combinations of the state variables. We assume that the displacement field belongs to the kernel of this form ($c$ linear form in Eq.\eqref{sec:ROMmethodo:subsec:formqs;subsubsec:continuous:eq:BCs}). In the general framework of the unidimensional problem, Neumann BCs on a given cable $\mathcal{C}_i$ are expressed as application of nodal forces $F_{i, j}$ applied on a set of discrete points $\{x^{\mathcal{C}_i}_j\}_{j=1}^{n_{\mathcal{C}_i}^{\rm 1d}}$. This translates into a jump $\llbracket . \rrbracket$ in the normal efforts at every point $x^{\mathcal{C}_i}_j$. In the end, the multi-modeling problem can be written as:

\begin{equation}
\label{sec:ROMmethodo:subsec:formqs;subsubsec:continuous:eq:varform}
\left\{
\begin{array}{rclrcl}
\mathcal{G}_\mu\left(\mathfrak{S}_\mu^{(k)}\right) &=& f^{(k)} \qquad &\text{on}&\quad  \ \Omega, \\
\mathfrak{S}_\mu^{(k)} &=& \StepQSOp^{(k)}\left(u^{(k)}_\mu, \ u^{(k-1)}_\mu,  \mathfrak{S}_\mu^{(k-1)}, \ H^{(k)} \right)\qquad &\text{on}&\quad  \ \Omega,  \\
\end{array}
\right.
\end{equation}

\noindent with BCs expressed as follows:

\begin{equation}
\label{sec:ROMmethodo:subsec:formqs;subsubsec:continuous:eq:BCs}
\resizebox{\columnwidth}{!}{$
\left\{\begin{array}{rcl}
\text{Dirichlet BCs}: &\quad & c(u_\mu^{(k)}) = 0\ \text{on} \ \Omega,\\ \\
\text{Neumann BCs}: &\quad &
\left\{
\begin{array}{rclrcl}
\left(\stressPCon\right)^{(k)} \cdot n &=& f_s^{(k)}\qquad &\text{on}&\quad  \ \Gamma_{\rm n}^{\rm c},  \\
\llbracket\stressPSte^{(k)}\rrbracket (x^{\mathcal{C}_i}_j) &=& F_{i, j}^{(k)} \qquad &\forall j\in \{1, ..., n_{\mathcal{C}_i}^{\rm 1d}\} &\quad \text{for} \quad \mathcal{C}_i, \ \forall i\in \{1, ..., n_{\mathcal{C}}\}, \\
\end{array}
\right.
\end{array}\right.$}
\end{equation}

Finally, the variational problem investigated in this contribution can be summarized as follows: eventually, the multi-modeling problem written in compact form in the Eq.\eqref{sec:ROMmethodo:subsec:formqs;subsubsec:continuous:eq:varform} to which the BCs are applied lead to the following variational problem:

\begin{equation}
\label{sec:ROMmethodo:subsec:formqs;subsubsec:continuous:eq:varform}
\resizebox{\columnwidth}{!}{$
\forall k\in \{1, ..., K\}, \ \text{Find} \ u_\mu^{(k)}\in \mathcal{X}_{\rm bc}\ \text{s.t.}
\left\{
\begin{array}{lcl}
\mathcal{R}_\mu\pr{u_\mu^{(k)}, \ u_\mu^{(k-1)}, \ \mathfrak{S}_\mu^{(k)}} =0,&\quad& \forall v\in \mathcal{X}_{\rm bc},\\
\mathfrak{S}_\mu^{(k)} = \StepQSOp^{(k)}\left(u^{(k)}_\mu, \ u^{(k-1)}_\mu,  \mathfrak{S}_\mu^{(k-1)}, \ H^{(k)} \right) &\text{on}& \Omega,\\
\end{array}
\right.$}
\end{equation}

\noindent where $\mathcal{X}_{\rm bc} \coloneqq \{v\in \mathcal{X}, c(v)=0\ \text{on} \ \Omega\}$. We denote:

\begin{equation*}
\resizebox{\columnwidth}{!}{$
\mathcal{R}_\mu\pr{u_\mu^{(k)}, \ u_\mu^{(k-1)}, \ \mathfrak{S}_\mu^{(k)}} = \mathcal{R}_\mu^{\mathfrak{S}}\pr{\StepQSOp^{(k)}\left(u^{(k)}_\mu, \ u^{(k-1)}_\mu,  \mathfrak{S}_\mu^{(k-1)}, \ H^{(k)}\right), \ v}, \ \text{and} \ \mathcal{R}_\mu^{\mathfrak{S}}\pr{\mathfrak{S}, \ v}=
\begin{bmatrix}
\mathcal{R}_\mu^{\sigma}\pr{\sigma_\mu^{(k)}, \ v} \\
\mathcal{R}_\mu^{\text{N}}\pr{\text{N}_\mu^{(k)}, \ v}
\end{bmatrix},$}
\end{equation*}

\noindent where we introduce the notations $\forall v\in [v^{\rm c}, \ v^{\rm s}]^\top$:

\begin{equation*}
\left\{\begin{array}{rcl}
\mathcal{R}_\mu^{\sigma}\pr{\sigma_\mu^{(k)}, \ v}&=&\displaystyle \int_{\Omega}\sigma_\mu^{(k)}:\varepsilon\pr{v^{\rm c}} \ d\Omega - \int_{\Omega}f_v\cdot v^{\rm c} \ d\Omega - \int_{\Gamma}f_s\cdot v^{\rm c} \ d\Gamma,  \\ \\
\mathcal{R}_\mu^{\text{N}}\pr{\text{N}_\mu^{(k)}, \ v}&=& \displaystyle \int_{\mathcal{C}} \text{N}_\mu^{(k)}:\partial_s v^{\rm s} \ ds - \int_{\mathcal{C}}f_v\cdot v^{\rm s} \ ds - \sum\limits_{i=1}^{n_{\mathcal{C}}}\sum\limits_{j=1}^{n_{\mathcal{C}_i}^{\rm 1d}} F_{i,j}^{(k)} v^{\rm s}(x_j^{\mathcal{C}_i}).
\end{array}\right.
\end{equation*}

\subsubsection{Finite element discretization}

We apply a problem discretization with a continuous Galerkin finite element (FE) method. Given the domain $\Omega$, we consider a HF mesh $\mathcal{T}^{\rm hf} =  \left\{\texttt{D}_i\right\}_{i=1}^{N_{\rm e}^{3d}} \cup \left\{\texttt{D}_i\right\}_{i=1}^{N_{\rm e}^{1d}}$ where $\texttt{D}_1^{1d},\ldots,\texttt{D}_{N_{\rm e}^{1d}}$ (resp. $\texttt{D}_1^{3d},\ldots,\texttt{D}_{N_{\rm e}^{3d}}$)  are the elements of the one dimensional-mesh, and $N_{\rm e}^{1d}$ (resp. $N_{\rm e}^{3d}$) denotes the number of elements in the one-dimensional (resp. three-dimensional) mesh. The $\rm hf$ subscript or superscript stands for HF discretization. In this framework, we denote by $\mathcal{X}^{\rm hf}$ the chosen finite element space to discretize the problem. Within this framework, we denote the displacement unknowns at nodes (primal) by $\FEdisP\in \mathbb{R}^{\mathcal{N}}$, where $\mathcal{N}=3(\mathcal{N}_{\rm no}^{\rm 3d}+\mathcal{N}_{\rm no}^{\rm 1d})$ is the dimension of the space $\mathcal{X}^{\rm hf}$. Furthermore, the generalized forces within the material are denoted by $\bm{\mathfrak{S}}_\mu=[\FEstressP, \ \fe{N}_\mu]^\top\in \mathbb{R}^{\mathcal{N}_g}$, since they are unknowns at quadrature points. For the record, the size of these vectors is $\mathcal{N}_g = \mathcal{N}_g^{\rm 3d}+\mathcal{N}_g^{1d} = 6\mathcal{N}_{\rm qd}^{\rm 3d}+\mathcal{N}_{\rm qd}^{1d}$, where $\mathcal{N}_{\rm qd}^{\rm 3d}$ stands for the number of quadrature weights used for the three-dimensional mesh and $\mathcal{N}_{\rm qd}^{1d}$ for the one-dimensional mesh. \\

We denote by $\{\mathbf{u}_\mu^{{\rm hf}, (k)}\}_{k=1}^K$ the FE approximation of the displacement (primal variable) given by the HF-model at all times, whereas $\{\bm{\mathfrak{S}}_\mu^{{\rm hf}, (k)}\}_{k=1}^K$ stand for the generalized force fields (stress or normal efforts). We state the FE discretization of the variational form Eq.\eqref{sec:ROMmethodo:subsec:formqs;subsubsec:continuous:eq:varform}, $\forall k\in \{1,..., K\}$, find $u^{{\rm hf}, (k)}_\mu\in \mathcal{X}_{\rm bc}^{\rm hf}$ s.t:

\begin{equation}
\label{sec:form:subsec:fe:eq:discritsystem}
\left\{
\begin{array}{rcl}
&\mathcal{R}^{\rm hf}_\mu\left(\mathbf{u}^{{\rm hf}, (k)}_\mu, \ \mathbf{u}^{{\rm hf}, (k-1)}_\mu, \ \bm{\mathfrak{S}}^{{\rm hf}, (k-1)}_\mu, \mathbf{v}\right)=0,& \qquad \forall v \in \mathcal{X}_{\rm bc}^{\rm hf},\\
&\bm{\mathfrak{S}}^{{\rm hf}, (k)}_\mu  = \mathcal{F}^{\rm hf}_\mu\left(\mathbf{u}^{{\rm hf}, (k)}_\mu, \mathbf{u}^{{\rm hf}, (k-1)}_\mu, \ \bm{\mathfrak{S}}^{{\rm hf}, (k-1)}_\mu, \ \fe{H}^{{\rm hf}, (k)}\right),&\\
\end{array}
\right.
\end{equation}

\noindent where $\mathcal{X}_{\rm bc}^{\rm hf} \coloneqq \left\{\mathbf{v}\in \mathcal{X}^{\rm hf}: \quad  \mathbf{B}\mathbf{v} = 0\right\}$ depicts the test space for displacements, and $\mathbf{B}\in \mathbb{R}^{\mathcal{N}_d\times\mathcal{N}}$ is the kinematic relationship matrix. $\mathcal{N}_d$ stands for the number of linear relations between degrees of freedom that we intend to enforce. Such a formulation on the BCs implies that the kinematic linear application depends neither on time nor on the parameter. Each line reflects a kinematic relationship between nodes of the overall mesh. Therefore, the said matrix includes not only the Dirichlet conditions applied to each physical domain, but also the kinematic relationships between the nodes of two distinct models (kinematic coupling). The operators $\mathcal{R}^{\rm hf}_\mu$ and $\mathcal{F}^{\rm hf}_\mu$ stands for the discrete counterparts of the continuous operators $\mathcal{R}_\mu$ and $\mathcal{F}_\mu$ introduced in Eq.\eqref{sec:form:subsec:fe:eq:discritsystem}. In practice, the FE code compute the HF-residuals as sums of elementary contributions, as follows $\forall v\in \mathcal{X}^{\rm hf}$:

\begin{footnotesize}
\begin{align*}
\label{sec:form:subsec:fe:eq:discritsystem:elementwise}
&\mathcal{R}^{\rm hf}_\mu\left(\mathbf{u}^{(k)}_\mu, \ \mathbf{u}^{(k-1)}_\mu, \ \bm{\mathfrak{S}}^{(k-1)}_\mu, \ \mathbf{v}\right)=\sum\limits_{q=1}^{N_e}\mathcal{R}^{\rm hf}_{\mu, q}\left(\mathbf{E}_q^{\rm no}\mathbf{u}^{(k)}_\mu, \ \mathbf{E}_q^{\rm no}\mathbf{u}^{(k-1)}_\mu, \ \mathbf{E}_q^{\rm qd}\bm{\mathfrak{S}}^{(k-1)}_\mu , \ \mathbf{E}_q^{\rm no}\mathbf{v}\right)\\
&=\underbrace{\sum\limits_{q=1}^{N_e^{\rm 3d}}\mathcal{R}^{\rm hf}_{\mu, q}\left(\mathbf{E}_q^{\rm no}\mathbf{u}^{(k)}_\mu, \ \mathbf{E}_q^{\rm no}\mathbf{u}^{(k-1)}_\mu, \ \mathbf{E}_q^{\rm qd, 3d}\bm{\sigma}^{(k-1)}_\mu , \ \mathbf{E}_q^{\rm no}\mathbf{v}\right)}_{\coloneqq\mathcal{R}^{\rm hf, 3d}_\mu\left(\mathbf{u}^{(k)}_\mu, \ \mathbf{u}^{(k-1)}_\mu, \ \bm{\sigma}^{(k-1)}_\mu, \ \mathbf{v}\right)} + \underbrace{\sum\limits_{q=1}^{N_e^{\rm 1d}}\mathcal{R}^{\rm hf}_{\mu, q}\left(\mathbf{E}_q^{\rm no}\mathbf{u}^{(k)}_\mu, \ \mathbf{E}_q^{\rm no}\mathbf{u}^{(k-1)}_\mu, \ \mathbf{E}_q^{\rm qd, 1d}\text{\textbf{N}}^{(k-1)}_\mu , \ \mathbf{E}_q^{\rm no}\mathbf{v}\right)}_{\coloneqq\mathcal{R}^{\rm hf, 1d}_\mu\left(\mathbf{u}^{(k)}_\mu, \ \mathbf{u}^{(k-1)}_\mu, \ \text{\textbf{N}}^{(k-1)}_\mu, \ \mathbf{v}\right)},
\end{align*}
\end{footnotesize}

\noindent where $\mathbf{E}_q^{\rm no}$ (resp. $\mathbf{E}_q^{\rm qd}$) is an elementary restriction operator on vectors at nodes (resp. quadrature points). For operators on vectors at quadrature points, we adopt the specific notation $\mathbf{E}_q^{\rm qd, 3d}$ (resp. $\mathbf{E}_q^{\rm no, 1d}$) for the case where the elements are three-dimensional (resp. one-dimensional). We emphasize that the assembly procedure can be split into two terms, a loop for concrete elements and a second loop for steel elements.\\

In our work, the boundary conditions are treated by dualization. Therefore, we introduce Lagrange multipliers, and we solve the following saddle-point problem: Find $(\fe{u}^{(k)}_\mu, \bm{\lambda}^{(k)}_\mu)\in \mathbb{R}^{\mathcal{N}}\times \mathbb{R}^{\mathcal{N}_d}$ s.t.:

$$
\left\{
\begin{array}{rcl}
\HFresDis\left(\fe{u}_\mu^{(k)}, \ \fe{u}_\mu^{(k-1)}, \ \bm{\mathfrak{S}}_\mu^{(k-1)}\right) + \mathbf{B}^\top \bm{\lambda}_\mu^{(k)}  &=&0, \\
\mathbf{B} \fe{u}_\mu^{(k)} &=& 0.
\end{array}
\right.
$$

This study expands upon a previously established framework for projection-based ROMs in nonlinear mechanics with internal variables, broadening its applicability to more intricate phenomena. Indeed, Eq.\eqref{sec:form:subsec:fe:eq:discritsystem} is formally similar to Eq.(9) of Reference  \cite{https://doi.org/10.1002/nme.7385}. However, the approach presented here allows for the partitioning of the domain into distinct regions, each characterized by a specific solid mechanics model. Additionally, our method broadens its scope to handle a diverse range of mechanical problems by incorporating auxiliary variables. These variables encompass fields such as temperature and water content, introducing influences on the constitutive equations, evolution equations of internal variables, and thereby, the mechanical state of the material.

\subsection{Projection-based model order reduction approach\label{sec:ROMmethodo:subsec:projROM}}

\begin{figure}[h!]
\begin{center}
\includegraphics[scale=0.85]{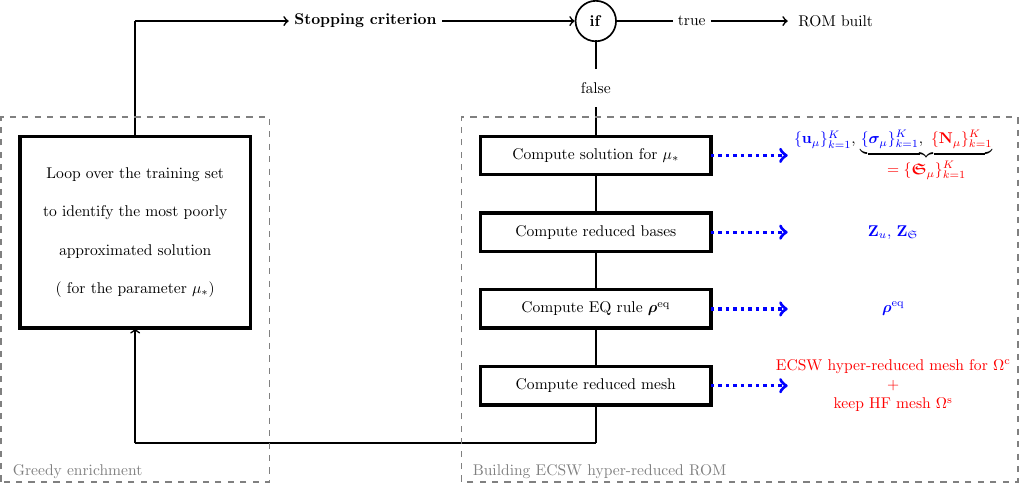} 
\end{center}
\caption{Key ideas of the greedy methodology and the ROM approach adopted for the mechanical simulation of prestressed concrete. The different components needed to define are provided in colors: \textcolor{blue}{blue} corresponds to components similar to the previous work, whereas \textcolor{red}{red} corresponds to the specific components within the multi-modeling framework for prestressed concrete \label{sec:method:subsec:projrom:fig:visustrat}}
\end{figure}

In this section, we discuss the pMOR procedure that is sketched in Figure \ref{sec:method:subsec:projrom:fig:visustrat}. As stated in the introduction, our method is an extension of the work \cite{https://doi.org/10.1002/nme.7385} to a more complex nonlinear mechanics problem with 3D-1D coupling. Our formulation relies on an offline-online computational decomposition (cf. Figure \ref{sec:intro:fig:context}): during the offline (training) stage, we solve the HF model for several parameter values to construct the reduced basis, the EQ rule and the associated mesh; during the online (prediction) stage, we call the surrogate model to approximate the solution. In section \ref{sec:ROMmethodo:subsec:projROM:subsubsec:srpb}, we consider the solution reproduction problem which addresses the task of reproducing the temporal trajectory for the same parameter value considered in the offline stage. In section \ref{sec:ROMmethodo:subsec:projROM:subsubsec:pROM}, we describe the extension to the parametric case.

\subsubsection{Solution Reproduction Problem\label{sec:ROMmethodo:subsec:projROM:subsubsec:srpb}}

We seek the reduced-order solution as a linear combination of modes:

\begin{equation*}
\widehat{\mathbf{u}}^{(k)}_\mu = \sum\limits_{n=1}^{N_u} \left(\widehat{\bm{\alpha}}^{(k)}_{u,\mu}\right)_n\bm{\zeta}_{u,n} = \mathbf{Z}_u\widehat{\bm{\alpha}}^{(k)}_{u,\mu},
\end{equation*}

\noindent where $\widehat{\bm{\alpha}}^{(k)}_{u,\mu}\in \mathbb{R}^{N_u}$ are referred to as generalized coordinates and $[\fe{Z}_u]_{.,n} = \bm{\zeta}_{u,n}$ are the displacement reduced basis vectors, which are build thanks to the POD approach. The main objective of this method is finding low dimensional approximations to the data, which preserve the essential information of a given high dimensional data set. More precisely, we resort to the method of snapshots \cite{sirovich1987turbulence} to build the displacement reduced basis. Given a discrete set of HF snapshots $\{\fe{v}_k\}_{k=1}^K$, a discrete scalar product $\left(\cdot, \cdot\right)$, and a tolerance $\varepsilon$, we define a Gramian matrix $\fe{C}\in \R^{K\times K}$, defined as $\fe{C}_{i,j} = (\fe{v}_i, \fe{v}_j)$. We need to solve a eigenvalue problem:

$$\fe{C}\bm{\varphi}_n = \lambda_n \bm{\varphi}_n, \quad \lambda_1 \geq \dotsc \geq \lambda_K, $$

\noindent to obtain the eigenpairs ($\lambda_n, \ \bm{\varphi}_n$). Thanks to the latter, we can compute POD modes:

$$\bm{\zeta}_{u,n} = \frac{1}{\sqrt{\lambda_n}}\sum\limits_{k=1}^K \left(\bm{\varphi}_n\right)_k \fe{v}_k.$$

\noindent The number of selected POD modes is chosen according to a energy-criterion on the spectrum, thanks to a user-defined tolerance $\varepsilon$:

\begin{equation*}
N_u = \min \left\{Q\in \mathbb{N}, \quad \sum\limits_{k=1}^Q\lambda_q \geq \left(1-\varepsilon^2\right)\sum\limits_{q=1}^K\lambda_q\right\}.
\end{equation*}

\noindent A use of the method of snapshots \cite{sirovich1987turbulence} for POD can thus be interpreted as a call to the following operator:

\begin{equation}
\label{eq:PODsnapOperator}
\fe{Z} = \text{POD}\left\{\{\fe{v}_k\}_{k=1}^K, \ \left(\cdot, \cdot\right), \ \varepsilon \right\}.
\end{equation}

The Galerkin ROM is obtained by projecting the discrete residual operator (onto the Eq.\eqref{sec:form:subsec:fe:eq:discritsystem}) onto the primal reduced basis. We first consider the situation without Lagrange multipliers for the boundary conditions:

\begin{equation}
\label{eq:ROBSolving}
\mathbf{Z}_u^\top\mathbf{R}^{\rm hf}_\mu\left(\widehat{\mathbf{u}}_\mu^{(k)}, \ \widehat{\mathbf{u}}_\mu^{(k-1)}, \ \widehat{\bm{\mathfrak{S}}}_\mu^{(k-1)}\right)=0.
\end{equation}

\noindent The nonlinearity of the operator results in a CPU bottleneck, since the assembly procedure scales with the cost of an HF computation. In order to circumvent this issue, we resort to an hyper-reduction approach, namely the element-wise EQ approach \cite{iollo2022adaptive}\cite{yano2019discontinuous}. The method samples a subset of the mesh elements over the entire computational domain in order to reduced the assembly costs in the online stage. With this approach, a residual operator $\mathcal{R}^{\rm eq}_\mu$ is generated and applied to the assembly procedure when the ROM solver is called. In the context of our multi-modeling problem, we choose to apply the hyper-reduction procedure only to three-dimensional terms, since these are nonlinear:

\begin{equation*}
\resizebox{\columnwidth}{!}{$
\mathcal{R}^{\rm eq}_\mu\left(\mathbf{u}^{(k)}_\mu, \ \mathbf{u}^{(k-1)}_\mu, \ \bm{\mathfrak{S}}^{(k-1)}_\mu, \ \mathbf{v}\right)=\underbrace{\mathcal{R}^{\rm eq, 3d}_\mu\left(\mathbf{u}^{(k)}_\mu, \ \mathbf{u}^{(k-1)}_\mu, \ \bm{\sigma}^{(k-1)}_\mu, \ \mathbf{v}\right)}_{\text{hyper-reduced}} +\underbrace{\mathcal{R}^{\rm hf, 1d}_\mu\left(\mathbf{u}^{(k)}_\mu, \ \mathbf{u}^{(k-1)}_\mu, \ \text{\textbf{N}}^{(k-1)}_\mu, \ \mathbf{v}\right)}_{\text{not hyper-reduced}}.$}
\end{equation*}

For the construction of the EQ rule, we rely on the Energy-Conserving Sampling and Weighting method (ECSW) developed in the Reference \cite{farhat2014dimensional}, whose quality has already been demonstrated for hyper-reduction of problems in solid mechanics. The ECSW approach consists in solving a non-negative  least-square problem to find a sparse approximation of the HF rule that is tailored to the integrals considered in \eqref{eq:ROBSolving}. Solving the optimization problem provides an EQ rule $\bm{\rho}^{\rm eq}\in \R^{N_e^{\rm 3d}}$, which defines the operator $\mathcal{R}^{\rm eq}_\mu$ from the HF operator as follows:

\begin{equation*}
\mathcal{R}^{\rm eq, 3d}_\mu\left(\mathbf{u}^{(k)}_\mu, \ \mathbf{u}^{(k-1)}_\mu, \ \bm{\sigma}^{(k-1)}_\mu, \ \mathbf{v}\right) = \sum\limits_{q=1}^{N_e^{\rm 3d}} \pr{\bm{\rho}^{\rm eq}}_q\mathcal{R}^{\rm hf}_{\mu, q}\left(\mathbf{E}_q^{\rm no}\mathbf{u}^{(k)}_\mu, \ \mathbf{E}_q^{\rm no}\mathbf{u}^{(k-1)}_\mu, \ \mathbf{E}_q^{\rm qd, 3d}\bm{\sigma}^{(k-1)}_\mu , \ \mathbf{E}_q^{\rm no}\mathbf{v}\right).
\end{equation*}

As already mentioned in the previous work by \citet{https://doi.org/10.1002/nme.7385}, the dualization of the BCs and the homogeneous BCs prevent us from considering any BCs in solving the problem, since the displacement modes satisfy the BCs (they belong to the $\fe{B}$ kernel). This highlights the fact that hyper-reduction of the three-dimensional domain, while preserving the one-dimensional part, has no impact on the application of the BCs. Information on the kinematic coupling between the steel and concrete nodes is already contained in the displacement modes. \\

Knowledge of the mechanical state of the material requires to know the stress field on the HF mesh. The latter is determined by integrating the constitutive equations at the quadrature points. However, the internal variables are only known at the sampled elements in the mesh. Hence, the stress field is only known at the reduced mesh level. To solve this problem, we build a reduced order basis for the generalized force $\bm{\mathfrak{S}}=[\bm{\sigma}, \text{\textbf{N}}]^\top$. Reconstruction of the generalized force field over the entire HF mesh is then performed using a Gappy-POD procedure \cite{everson1995karhunen}. Unlike displacement vectors, the components of generalized force vectors on one-dimensional and three-dimensional discrete points do not have the same physical dimension \cite{taddei2021discretize}\cite{parish2023impact}. Therefore, we define the scalar product on generalized force vectors:

\begin{equation*}
\left(\bm{\mathfrak{S}}_1, \bm{\mathfrak{S}}_2\right) = \left(\begin{bmatrix} \bm{\sigma}_1 \\ \fe{N}_1 \end{bmatrix}, \ \begin{bmatrix} \bm{\sigma}_2 \\ \fe{N}_2 \end{bmatrix}\right)_{[\sigma, N]} = \frac{1}{\lambda_1^{\sigma}} \left(\bm{\sigma}_1, \ \bm{\sigma}_2\right)_2 + \frac{1}{\lambda_1^{\text{N}}} \left(\fe{N}_1, \ \fe{N}_2\right)_2,
\end{equation*}

\noindent where $\lambda_1^{\sigma}$ (resp. $\lambda_1^{\text{N}}$) is the largest eigenvalue in the sense of the scalar product $\ell_2$ for the stress vectors (normal forces). In summary, in addition to the EQ rule $\bm{\rho}^{\rm eq}$ (and the associated reduced mesh), the ROM is made up of two reduced bases, defined thanks to the POD operator detailed in Eq.\eqref{eq:PODsnapOperator} as follows:

\begin{equation*}
\fe{Z}_u = \text{POD}\left\{\{\fe{u}^{{\rm hf}, (k)}_\mu\}_{k=1}^K, \ \left(\cdot, \cdot\right)_2, \ \varepsilon_u \right\}, \quad \text{and} \quad \fe{Z}_{\bm{\mathfrak{S}}} = \text{POD}\left\{\{\bm{\mathfrak{S}}_\mu^{{\rm hf}, (k)}\}_{k=1}^K, \ \left(\cdot, \cdot\right)_{[\sigma, N]}, \ \varepsilon_{\mathfrak{S}} \right\}.
\end{equation*}

Both for the displacements and for the generalized forces, we opted for a scalar product $\ell_2$ on the discrete snapshots. From a variational perspective, it would have been more suitable to work with an H1 scalar product. However, from an algorithm application point of view, extracting such matrices in this context can be fairly challenging. Our choice is typical for numerical applications on real-world applications. Furthermore, $\ell_2$ compression delivers high-quality numerical results.

\subsubsection{Parametric problem\label{sec:ROMmethodo:subsec:projROM:subsubsec:pROM}}

In order to provide a reliable ROM on a set of parameters, we build the surrogate model using a POD-Greedy approach. This iterative procedure is designed to enrich the reduced model (i.e. the reduced bases and the reduced mesh) by computing at each iteration the HF solution least well approximated by the ROM. The worst-approximated solution is estimated by exploring a test set $\Theta_{\rm train}$, defined as a discrete approximation of $\mathcal{P}$. In our case, we chose to rely on a strong-greedy approach: we compare the approximation errors (error between HF solution and reduced solution) over the whole test set, to identify the parameter for which this error is maximal. This parameter is then used to further enhance the ROM. Strong-greedy approach is not optimal from the standpoint of the computational cost of building the reduced model. Indeed, the estimation of the poorest approximated solution requires knowledge of the HF solutions on a given discrete training set. For a more efficient greedy approach in terms of computational cost, weak-greedy methods would be more appropriate, along with the introduction of an appropriate error indicator. This remains a limitation to be borne in mind, particularly in the context of increasing the dimensionality of the parameter space. Nevertheless, this work constitutes a proof of concept of the feasibility of a greedy approach for three-dimensional THM calculations on prestressed concrete. The numerical optimization of the process, with the development of error indicators adapted to these problems and to industrial-grade HF codes, is a focus for forthcoming research. \\

The switch of the methodology to the parametric case requires the adaptation of two parts of the algorithm: the construction of the reduced bases and the computation of the EQ. Two constructions of the reduced bases are explored in this paper. A first approach consists in performing a new POD on the set of computed HF snapshots. A second approach involves an incremental approach, known in the literature as H-POD \cite{haasdonk2008reduced}. The latter has the advantage of providing a hierarchical basis achieved by concatenating the previous basis with one obtained with new snapshots:

\begin{equation*}
\fe{Z} = \left[\fe{Z}, \ \fe{Z}_{\rm proj}\right], \quad \fe{Z}_{\rm proj} = \text{POD}\left\{ \left\{\Pi_{\fe{Z}^\bot}\fe{v}_k \right\}_{k=1}^K
, \ \pr{\cdot, \cdot}, \ \varepsilon\right\},
\end{equation*}

\noindent where $\Pi_{\fe{Z}^\bot} \ : \ \mathcal{X}^{\rm hf} \rightarrow \mathcal{Z}$ is the orthogonal projection operator $\mathcal{Z}\subset \mathcal{X}^{\rm hf}$ using the $(\cdot, \cdot)$ scalar product. We rely on the regularization approaches given in Reference \cite{haasdonk2008reduced} to compute the new number of modes. Before concatenating the two bases, a criterion is added such that only the basis vectors that effectively reduce the projection error are added to the reduced order basis. The different steps of the adaptive algorithm are summarized in Algorithm \ref{alg:podgreedy}.\\

\begin{algorithm}[h]
\caption{strong POD-Greedy algorithm}\label{alg:podgreedy}
\begin{algorithmic}[1]
\Require $\Theta_{\rm train}=\{\mu_i\}_{i}^{n_{\rm train}}$, $\varepsilon_{u}$, $\varepsilon_{\mathfrak{S}}$
\State $\mathcal{Z}_{N_u} = \mathcal{Z}_{N_\sigma} = \emptyset$, $\mu_* = \overline{\mu}$, $\Theta_{*}=\{\mu_*\}$.
\While{ Stop Criterium }
\State Compute $\{\mathbf{u}^{{\rm hf}, (k)}_{\mu^*}\}_{k=1}^K$,  $\{\bm{\mathfrak{S}}^{{\rm hf}, (k)}_{\mu^*}\}_{k=1}^K$\Comment{Call of \textsf{code$\_$aster}}
\State Compute primal reduced basis $\mathbf{Z}_u$ 
\State Compute $\bm{\rho}^{\rm eq}$ knowing $\{\zeta_{u,n} \}_{n=1}^{N_u}$ and $\{\bm{\mathfrak{S}}^{{\rm hf}, (k)}_{\mu}\}_{k\in \{1, .., K\}, \mu\in \Theta_{\rm *}}$ 
\State Compute the reduced mesh $\mathcal{T}^{\rm red}$
\State Compute dual reduced basis $\mathbf{Z}_{\mathfrak{S}}$
\For{$\mu \in \Theta_{\rm train}$}
\State Solve the ROM for $\mu$ and compute $E_{\mu}^{\rm app, avg}$\Comment{See definition of $E_{\mu}^{\rm app, avg}$ in Eq.\eqref{sec:numres_srpb:approxError}}
\EndFor
\State $\mu^* = \arg  \max\limits_{\mu\in \Theta_{\rm train}}E_{\mu}^{\rm app, avg}$
\State $\Theta_{\rm *}= \Theta_{\rm *}\cup \{\mu_*\}$
\EndWhile
\end{algorithmic}
\end{algorithm}

\section{Thermo-Hydro-Mechanical (THM) modeling of large concrete structures\label{sec:THMmodel}}
\subsection{Weak-coupling strategy for the THM numerical model}

In this section, we introduce the mathematical model designed to simulate the behavior of prestressed concrete. We consider models that account for the evolution of large concrete structures over their lifetime, which consists mainly of two stages: young age and long-term evolution. The young age refers to a stage during which the chemical and physical properties of concrete are changing at a fast rate, as it sets and hardens. Long-term phase represents the evolution of hardened concrete under operating conditions (taking into account thermo-hydric loadings) and mechanical loadings. Within the framework of the FE models employed in practice, we only consider long-term evolution. The behavior of heterogeneous and porous concrete is governed by numerous and complex physicochemical phenomena. Such material behavior requires a THM modeling strategy: the behavior of the material is based on knowledge of the temperature ($T$), the water content ($C_w$) in the concrete and the mechanical fields, in a framework where all these phenomena are coupled together. Since we are interested in modeling the whole ageing of the concrete structure, our THM model should encompass the various physical processes which induce deformations within the concrete: shrinkage, dessication and creep. \\

\begin{figure}[h!]
\begin{center}
\includegraphics[scale=0.95]{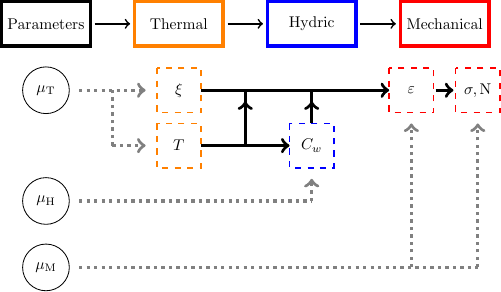} 
\end{center}
\caption{Weakly-coupled chained THM approach for large prestressed concrete structures\label{sec:THMmodel:subsec:WC:fig:WC}. Each step provides different fields of interest: at the end of the thermal calculation, we get the temperature field ($T$) and the degree of hydration of the concrete ($\xi$; which will be always analytically given in our simulations); at the end of the hydric calculation, we get the water content of concrete ($C_w$); at the end of the mechanical calculation, we get the displacement fields in the steel cables and in the concrete, the associated deformations ($\varepsilon = [\varepsilon^{\rm c}, \varepsilon^{\rm s}]$), the stresses in the concrete ($\sigma$) and the normal forces in the cables ($\text{N}$).}
\end{figure}

\begin{table}[h!]
\begin{center}
\begin{tabular}{ccc}
\hline
\textbf{Notation} & \textbf{Physical quantity} & \textbf{Unit} \\ \hline
$T$ & Temperature &  \si{K} \\ 
$\xi$ & Hydration degree &  \si{-} \\ \hdashline
$h$ & Ambient relative humidity (${\rm RH}$) & \si{-} \\
$C_w$ & Water content of concrete & \si{-} \\ \hdashline
$\sigma$ & Stress in the concrete & \si{Pa} \\
$\varepsilon^{\rm c} = \nabla_s u^{\rm c}$ & Deformations in the concrete & \si{-} \\
$N$ & Normal efforts in the prestressing cables & \si{N} \\
$\varepsilon^{\rm s} = \partial_s u^{\rm s}$ & Deformations in the prestressing cables & \si{-} \\
\end{tabular}
\end{center}
\caption{Fields of interest in the overall THM model for large prestressed concrete structures\label{sec:THMmodel:subsec:WC:table_fields}}
\end{table}
In our framework, we adopt a weakly-coupled approach. This assumption implies that the computation is carried out in a chained manner. First, a thermal calculation is performed, followed by a hydric calculation (water diffusion in the concrete). Once all the thermal and hydric fields are known, a mechanical calculation is conducted. Each calculation step yields fields of interest which: first describe the state of the material; second, can be reused for subsequent calculation steps. The Table \ref{sec:THMmodel:subsec:WC:table_fields} details the ouputs for the entire THM calculation. The different steps in the process are summarized in Figure \ref{sec:THMmodel:subsec:WC:fig:WC}. Such a formulation of the problem is founded on several assumptions. To begin with, the influence of the mechanical response on the thermal and water fields is neglected \cite{jason2007hydraulic}\cite{bouhjiti2018accounting}. Furthermore, it is assumed that the hydric response has no influence on the thermal fields. Weak-coupling approaches have demonstrated their effectiveness in modeling prestressed concrete structures, both for a RSV \cite{jason2007hydraulic}\cite{bouhjiti2018accounting} and for a full-scale model \cite{asali2016numerical}. 

\subsection{THM constitutive equations\label{sec:THMmodel:subsec:consEq}}

As stated above, we describe in the following section the set of equations that make up the THM problem under study. Prestressed concrete behavior modeling requires a multi-modeling approach: a three-dimensional nonlinear rheological model is used for concrete; and prestressing cables are described by a one-dimensional linear thermo-elastic behavior. As mentioned above, the rheological behavior of concrete is coupled with hydric and thermal phenomena. Thermal-hydric resolutions are thus solved on the concrete domain ($\Omega^{\rm c}$), while mechanical calculations are solved on both domains ($\Omega^{\rm c}$ and $\Omega^{\rm s}$).

\subsubsection{Modeling of the thermal and the hydric behavior of the concrete\label{sec:THMmodel:subsec:consEq:subsubsec:THEq}}

First, we introduce the set of equations employed for the first two stages of the chained calculation: the thermal calculation and the hydric calculation. This decision is motivated by the fact that this calculation is the starting point for the mechanical calculation, to which we apply our model reduction methodology (section \ref{sec:ROMmethodo}). The temperature evolution is modeled by the heat equation \cite{fourier1822theorie}:

\begin{equation}
\label{sec:THMmodel:subsec:consEq:subsubsec:THEq:eq:Teq}
 \rho_c c_p^p \frac{\partial T}{\partial t} = \nabla \cdot \left(\lambda_c \nabla T\right), \quad \text{on} \ \Omega_{\rm c},
\end{equation}

\noindent where $\rho_c$ is the density of the concrete, $c_p^p$ heat capacity of hardened concrete and $\lambda_c$ thermal conductivity of hardened concrete. Dirichlet conditions are applied in our context (see details for the numerical test case in section \ref{sec:THMmodel:subsec:RSV}).\\

Since we only consider liquid water diffusion \cite{granger1995comportement}, moisture transfer is modeled by a single nonlinear diffusion on $C_w$ (see Eq \eqref{sec:THMmodel:subsec:consEq:subsubsec:THEq:eq:hydric_eq}), which denotes the water content of the concrete. The diffusion equation depends on $D_w$, which is a phenomenological diffusion coefficient, and is assumed to follow Arrhenius' law\cite{bavzant1972nonlinear}. In summary, the nonlinear diffusion equation of the water content can be summed up as follows:

\begin{subequations}
\begin{empheq} [left=\empheqlbrace] {align}
\frac{\partial C_w}{\partial t} &= \nabla \cdot \left[D_w\left(C_w, \ T\right) \nabla C_w\right], \quad \text{on} \ \Omega_{\rm c,}\label{sec:THMmodel:subsec:consEq:subsubsec:THEq:eq:hydric_eq} \\
D_w\left(C_w, \ T\right) &=  D_{w,0} \left(C_w\right) \frac{T}{T_w^0} \exp\left[-\frac{U_w}{R}\left(\frac{1}{T} - \frac{1}{T_w^0}\right)\right],\label{sec:THMmodel:subsec:consEq:subsubsec:THEq:eq:hydric_eq_2}\\
D_{w,0}\left(C_w\right) &= A \exp \left(B C_w\right),\label{sec:THMmodel:subsec:consEq:subsubsec:THEq:eq:hydric_eq_3}
\end{empheq}
\end{subequations}

\noindent where $U_w$ is the activation energy of drying, $R$ the ideal gas constant and $D_{w,0}\left(C_w\right)$ is the diffusion coefficient at a reference temerature $T_w^0$. The latter is assumed to follow a model defined by \citet{mensi1988sechage}, which depends on two model parameters $A$, and $B$. \\

At the scale of large concrete structures, measurements of ambient conditions cannot be made in terms of the water content of the concrete, and are thus conducted in relative humidity \cite{boucher2016analyse}. Relative humidity (RH) is defined as the ratio of vapor pressure to saturation vapor pressure for a given temperature. For the sake of consistency and use of collected data, the boundary conditions are formulated in terms of RH.  From an experimental point of view, the drying or wetting cycles are assumed to affect only the concrete skin. This assumption enables to draw a link between the water concentration in the concrete and the relative humidity. For a given temperature, these two quantities are related by a bijective function called the sorption-desorption function:

\begin{equation}
\label{sec:THMmodel:subsec:consEq:subsubsec:THEq:eq:sorptionfunc}
C_w = f_d\left(h\right).
\end{equation}

\begin{figure}[h!]
\begin{center}
\begin{subfigure}[b]{0.45\textwidth}
\begin{center}
\includegraphics[scale=1]{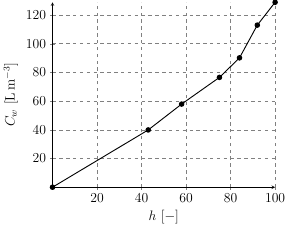} 
\end{center}
\caption{Soprtion-desorption function used for the THM problem}
\end{subfigure}
\qquad
\begin{subfigure}[b]{0.45\textwidth}
\begin{adjustbox}{width=\columnwidth,center}
\begin{tabular}{ccccccccc}
\hline $C_w$ [\si{L.m^{-3}}] & 0 & 39.0 & 57.9 & 76.5 & 90.1 & 112.9 & 128.8 \\
\hline $h$ [\si{-}] & 0 & 43 & 58 & 75 & 84 & 92 & 100 
\end{tabular}
\end{adjustbox}
\caption{\normalsize Summary}
\end{subfigure}
\end{center}
\caption{Definition of the sorption-desorption function $f_d$ (defined in Eq \eqref{sec:THMmodel:subsec:consEq:subsubsec:THEq:eq:sorptionfunc}). The table shows the point values given to define the function. The function is computed by linear interpolation between those points. The reference configuration corresponds to $h=100$, which is the initial RH value in the wall.\label{sec:THMmodel:subsec:consEq:subsubsec:THEq:hydric:fig:soprtionFunc}}
\end{figure}

Within the framework of these constitutive laws, the sorption-desorption function may be defined either analytically with hyper-parameters \cite{bouhjiti2018accounting}\cite{{van1980closed}}, or empirically by defining a function. In our case, we define a sorption-desorption function as shown in Figure \ref{sec:THMmodel:subsec:consEq:subsubsec:THEq:hydric:fig:soprtionFunc}. This curve is drawn from experimentally acquired points (without interpolation).

\noindent As previsouly mentioned, the BC of the water diffusion problem are stated in terms of RH when using real life data. All the parameters related to the thermal and hydric aspects of the model are summarized and detailed in Table \ref{sec:THMmodel:subsec:consEq:subsubsec:THEq:table:parameters}.

\begin{table}[h!]
\begin{center}
\begin{adjustbox}{width=\columnwidth,center}
\begin{tabular}{cccc}
\hline
\textbf{Calculation step }& \textbf{Notation} & \textbf{Physical quantity or parameter} & \textbf{Unit} \\ \hline
&$\rho_c$ & Density & \si{kg.m^{-3}} \\ 
Thermal ($\mu_{\text{T}}$)  &$c_c^p$ & Heat capacity of hardened concrete & \si{kJ.kg.K^{-1}}  \\ 
&$\lambda_c$ & Thermal conductivity of hardened concrete & \si{W.m^{-1}.K^{-1}} \\ \hdashline
&$D_w$ & Phenomenological diffusion coefficient &  \\
&$f_d$ & Sorption-desorption function & \\ 
&$T_{w}^0$ & Reference temperature & \si{K} \\
&$D_{w,0}$ & Diffusion coefficient at a reference temperature $T_{w}^0$ &  \\
&$U_w$ & Activation energy of drying & \si{kJ.mol^{-1}} \\
&$R$ & Ideal gas constant & \si{kJ.mol^{-1}.K^{-1}} \\ 
Hydric ($\mu_{\text{H}}$) &$A$ & Model parameter for Mensi's law & \si{10^{-15}m^2.s^{-1}} \\
&$B$ & Model parameter for Mensi's law & \si{-} \\
 \hline
\end{tabular}
\end{adjustbox}
\end{center}
\caption{Summary of parameters and physical quantities at stake in the modeling of the thermal (see Eq.\eqref{sec:THMmodel:subsec:consEq:subsubsec:THEq:eq:Teq} and the hydric (see Eq.\eqref{sec:THMmodel:subsec:consEq:subsubsec:THEq:eq:hydric_eq}-\eqref{sec:THMmodel:subsec:consEq:subsubsec:THEq:eq:hydric_eq_2}-\eqref{sec:THMmodel:subsec:consEq:subsubsec:THEq:eq:hydric_eq_3}) behavior \label{sec:THMmodel:subsec:consEq:subsubsec:THEq:table:parameters}}
\end{table}

\subsubsection{Modeling of the mechanical behavior of the concrete}

\begin{figure}[h!]
\begin{center}
\centering
\begin{subfigure}[b]{0.45\textwidth}
\includegraphics[scale=0.9]{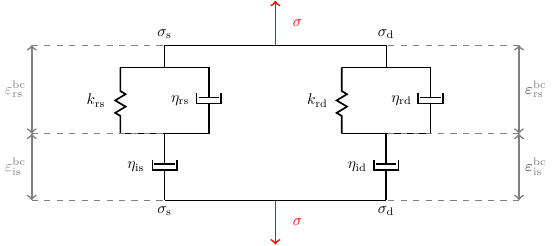} 
\caption{Burger rheological model for the basic creep}
\end{subfigure}
\qquad
\centering
\begin{subfigure}[b]{0.45\textwidth}
\begin{adjustbox}{width=\columnwidth,center}
\begin{tabular}{ccc}
\hline
\textbf{Notation} & \textbf{Physical quantity or parameter} & \textbf{Unit} \\ \hline
$E_{\rm c}$ & Young's modulus (concrete) & \si{Pa} \\ 
$\nu_{\rm c}$ & Poisson's ratio (concrete) & \si{-}  \\ 
$\alpha_{\rm th, c}$ & Thermal dilation coefficient (concrete) & \si{K^{-1}} \\ 
$\alpha_{\rm dc}$ & Dessication shrinkage coefficient & \si{-} \\
$\beta_{\rm endo}$ & Autogeneous shrinkage coefficient & \si{-} \\ 
$\nu_{\rm bc}$ & Basic creep Poisson's ratio & \si{-} \\
$k_{\rm rd}$ & Reversible deviatoric basic stiffness & \si{Pa} \\
$\eta_{\rm rd}$ & Reversible deviatoric basic viscosity & \si{Pa.s} \\
$\eta_{\rm id}$ & Irreversible deviatoric basic viscosity & \si{Pa.s} \\
$U_{\rm bc}$ & Basic creep activation energy & \si{kJ.mol^{-1}} \\ 
$T_{\rm bc}^0$ & Basic creep reference temperature & \si{^\circ C} \\ 
$\kappa$ & Basic creep consolidation parameter & \si{-} \\
$\eta_{\rm dc}$ & Desiccation creep parameter & \si{Pa.s} \\
\end{tabular}
\end{adjustbox}
\caption{\normalsize Summary of the parameters for the mechanical model}
\end{subfigure}
\end{center}
\caption{Parameters for the three-dimensional mechanical model (concrete)}
\end{figure}

In this section, we detail the governing equations for the mechanical behavior of concrete. Since, we consider small-displacement small-strain mechanical problems, the total strain is decomposed as the sum of several contributions:

\begin{equation*}
\varepsilon = \varepsilon^{\rm el} + \varepsilon^{\rm th} +  \varepsilon^{\rm ds} + \varepsilon^{\rm bc} + \varepsilon^{\rm dc},
\end{equation*}

\noindent where $\varepsilon^{\rm el}$ is the elastic strain tensor, $\varepsilon^{\rm th}$ the thermal strain tensor, $\varepsilon^{\rm ds}$ the dessication shrinkage strain tensor, $\varepsilon^{\rm bc}$ the basic creep strain tensor, $\varepsilon^{\rm dc}$ the dessication creep strain tensor and $\varepsilon^{\rm en}$ the autogenous shrinkage. We explain in the following section the evolution and constitutive equation that help expressing the different strain tensors.\\

According to experimental observations, the variation of thermal strain $\varepsilon^{\rm th}$ proportionnal to temperature variations (see Eq \eqref{sec:THMmodel:subsec:consEq:subsubsec:MEq:eq:th}). The proportionality coefficient $\alpha_{\rm th, c}$ is referred to as the thermal dilation coefficient of concrete and is assumed to be constant when focusing on the long-terme phase. Similar experimental observations suggest a linear dependency between the variations of the dessication shrinkage strains $\varepsilon^{\rm ds}$ and the water content of the concrete $C_w$ (see Eq \eqref{sec:THMmodel:subsec:consEq:subsubsec:MEq:eq:ds}), which is expressed thanks to dessication shrinkage coefficient ($\alpha_{\rm ds }$). We assume that we have the same kind of relationship between the autogenous shrinkage $\varepsilon^{\rm en}$ and the hydratation degree $ \xi$, expressed thanks to a $\beta_{\rm endo}$ coefficient.

\begin{subequations}
\begin{empheq} [left=\empheqlbrace] {align}
\dot{\varepsilon}^{\rm th} &= \alpha_{\rm th, c} \frac{\partial T}{\partial t} \text{I},\label{sec:THMmodel:subsec:consEq:subsubsec:MEq:eq:th} \\
\dot{\varepsilon}^{\rm ds} &=  \alpha_{\rm dc } \frac{\partial C_w}{\partial t} \text{I},\label{sec:THMmodel:subsec:consEq:subsubsec:MEq:eq:ds}\\
\dot{\varepsilon}^{\rm en} &=  \beta_{\rm endo} \frac{\partial \xi}{\partial t} \text{I}.\label{sec:THMmodel:subsec:consEq:subsubsec:MEq:eq:endo}
\end{empheq}
\end{subequations}

The model selected for the creep deformations is the Burger model developed by \citet{foucault2012new}. This choice is motivated by several experimental validations and is well-suited for creep investigations on the considered structures, as confirmed by the work of \citet{bouhjiti2018accounting}. We assume that the creep is a phenomenon involving a decoupling of a spherical part and a deviatoric part. We decompose the Cauchy stress tensor ($\sigma$) as the sum of a spherical part ($\sigma_{\rm s}$) and deviatoric part ($\sigma_{\rm d}$):

\begin{equation*}
\sigma = \sigma_{\rm s} \text{I} + \sigma_{\rm d},
\end{equation*}

\noindent where $\sigma_{\rm s} = \text{Tr}(\sigma)/3$, and $\text{I}$ is the identity tensor. The Burger creep model is built on a decomposition into a reversible and an irreversible part, where we split each tensor into its spherical and deviatoric part:

\begin{equation*}
\left\{
\begin{array}{rcl}
\varepsilon &=& \varepsilon^{\rm bc}_{\rm i} + \varepsilon^{\rm bc}_{\rm r},\\
\varepsilon^{\rm bc}_{\rm i} &=& \varepsilon^{\rm bc}_{\rm rs}\text{I} + \varepsilon^{\rm bc}_{\rm rd}, \\
\varepsilon^{\rm bc}_r &=& \varepsilon^{\rm bc}_{\rm is}\text{I} + \varepsilon^{\rm bc}_{\rm id}.\\
\end{array}
\right.
\end{equation*}

Each part (deviatoric and spherical) is described by a Burger-type model. For each chain, the reversible basic creep strains are modeled by  a Kelvin-Voigt rheological elements, whereas the irreversible basic creep strains are modeled by Maxwell elements. The Kelvin-Voigt model (see Eq. \eqref{sec:THMmodel:subsec:consEq:subsubsec:MEq:eq:spherical_rev}) used for the reversible reversible spherical basic creep is expressed thanks to the stiffness (resp. viscosity) $k_{\rm rs}$ (resp. $\eta_{\rm rs}$), while the irreversible spherical basic creep viscosity $\eta_{\rm is}$ is given by a nonlinear relationship, expressed thanks to a consolidation parameter $\kappa$ (see Eq.\eqref{sec:THMmodel:subsec:consEq:subsubsec:MEq:eq:spherical_irr}).

\begin{subequations}
\begin{empheq} [left=\empheqlbrace] {align}
h \sigma_{\rm s} &= k_{\rm rs} \varepsilon^{\rm bc}_{\rm rs} + \eta_{\rm rs} \dot{\varepsilon}^{\rm bc}_{\rm rs}, \label{sec:THMmodel:subsec:consEq:subsubsec:MEq:eq:spherical_rev} \\
h\sigma_{\rm s} &= \underbrace{\eta_{\rm is}^0 \exp \left(\frac{\norm{\varepsilon_i^{\rm bc}}_m}{\kappa}\right)}_{\coloneqq \eta_{\rm is}}\dot{\varepsilon}^{\rm bc}_{\rm is}, \ \text{where} \ \norm{\varepsilon_i^{\rm bc}}_m = \max\limits_{\tau \in \left[0, t\right]} \sqrt{\varepsilon_i^{\rm bc}\left(\tau\right):\varepsilon_i^{\rm bc}\left(\tau\right)}, \quad \forall t\geq 0,\label{sec:THMmodel:subsec:consEq:subsubsec:MEq:eq:spherical_irr} 
\end{empheq}
\end{subequations}

\noindent The deviatoric part is expressed using a similar set of tensor equations (the spherical part being a set of scalar equations). The aforementioned model accounts for thermo-activation of basic creep. To this end, stiffness and viscosity parameters expressions follow an Arrhenius' law:

\begin{equation*}
\kappa \left(T\right) = \kappa_0 \exp\left[-\frac{U_{\rm bc}}{R}\left(\frac{1}{T} - \frac{1}{T^0_{bc}}\right)\right],
\end{equation*}

\noindent where $k_{\rm rs}^0$ is the reversible spherical creep stiffness at a reference temperature $T_{\rm bc}^0$ and $U_{\rm bc}$ the activation energy of basic creep. Finally, the equivalence of spherical and deviatoric chains enables to restrict the number of model parameters, by assuming a constant creep Poisson ratio $\nu_{\rm bc}$, given by the following relation:

\begin{equation*}
\frac{k_{\rm rs}}{k_{\rm rd}} = \frac{\eta_{\rm rs}}{\eta_{\rm rd}} = \frac{\eta_{\rm rs}^0}{\eta_{\rm rd}^0} = \frac{1 + \nu_{\rm bc}}{1 - 2\nu_{\rm bc}}.
\end{equation*}

In order to model the dessication creep strain, we consider the following equation, founded on the work of \citet{bazant1985concrete}:
\begin{equation*}
\dot{\varepsilon}^{\rm bc} = \frac{1}{\eta_{\rm dc}}\left|\dot{h}\right| \sigma,
\end{equation*}

\noindent where $\eta_{\rm dc}$ is a material parameter (\si{Pa.s}).

\subsubsection{Modeling of the coupling between concrete and prestressing cables}

Within large prestressed concrete structures, the aim of the steel cables embedded within the concrete is to apply permanent compressive stresses to the concrete in such a way as to compensate for the tensile forces that are to be applied to the structure. This technique generates favorable internal forces in the concrete. The installation of a prestressed concrete structure requires the tensioning of the cables in the concrete. The tension profile along a cable is designed in our case to comply with an official standard (BPEL 91 regulation). Given physical parameters (initial tension), a tension profile is computed along the length of the cable as a function of the curvilinear abscissa. To be more precise, the coupling between cables and concrete can be decomposed into three main stages: before, during and after prestressing. In the case of the above structures, the concrete is first poured around sheaths and begins to dry. Cables are then inserted into these ducts and prestressed in order to comply with civil engineering standards. At last, cement is poured in the ducts, and the life of the structure can continue with a kinematic coupling between the concrete and the steel cables. In the numerical model studied here, the one-dimensional mesh (modeling the steel cables) is immersed within the three-dimensional mesh. This means that the cables "cross" the concrete cells. A kinematic linkage is performed in order to connect the concrete nodes and the steel nodes.  Since the coupling is assumed to be perfect (no slip between the tendons and the cement), coincident points in each material are assumed to have the same displacement. Instantaneous prestressing losses due to anchor recoil and friction are not taken into account at the scale of the considered RSV. \\

The cables are modeled by bars, which means that we resort to a one-dimensional approach where only the tension-compression forces are considered. In this framework, the structure is described at each instant by a curve representing its mean line. Consequently, only the normal efforts appear (efforts defined along the tangent vector to the beam section) in the variational formulation of the problem. Two sets of equilibrium equations appear in the studied case: during the prestressing step (namely between the times $t^{\rm init, p}$ and $t^{\rm end, p}$) and after the prestressing step (namely until the end of the study $t_{\rm f}$):

\begin{equation}
\label{eq:loadings_for_prestressed_step:anal}
\left\{
\begin{array}{rcl}
\partial_s N \left(s, \ t\right) &=& f_{{\rm s}}, \quad \forall t \in \left[t^{\rm init, p}\ , t^{\rm end, p}\right], \quad \text{and}\quad \llbracket N \rrbracket\pr{x_i^{\rm no, 1d}} = - \frac{t^{(k)} - t^{\text{init,p}}}{t^{\text{end,p}} - t^{\text{init,p}}}F_i, \\
\partial_s N \left(s, \ t\right) &=& f_{{\rm s}}, \quad \forall t \in \left[t^{\rm end, p}, t_{\rm f}\right], \\
\end{array}\right.
\end{equation}

\noindent where $\{x_i^{\rm no, 1d}\}_{i=1}^{\mathcal{N}^{s}}$ are the nodes of the one-dimensional mesh and $F_i$ are the nodal forces prescribed in order to respect the BPEL regulation used. We consider a linear thermo-elastic constitutive equation for the steel cables. Thus, the normal efforts in the cables ($\text{N}$) are linked to the uniaxial strains ($\varepsilon^{\rm s}$) in the cables:

\begin{equation*}
\text{N} = E_s S_s \left(\varepsilon^{\rm s} - \alpha_{\rm th, s} \Delta T\right),
\end{equation*}
\noindent where $E_s$ the Young's modulus, $\alpha_{\rm th, s}$ the thermal dilation coefficient, $S_s$ the section of the prestressing cables and $\Delta T$ is  the temperature rise in the beam.\\

\begin{figure}[h!]
\begin{center}
\begin{tabular}{ccc}
\hline
\textbf{Notation} & \textbf{Physical quantity or parameter} & \textbf{Unit} \\ \hline
$E_{\rm s}$ & Young's modulus (steel) & \si{Pa} \\ 
$\nu_{\rm s}$ & Poisson's ratio (steel) & \si{-}  \\ 
$\rho_{\rm s}$ & Density (steel) & \si{kg.m^3} \\
$\alpha_{\rm th, s}$ & Thermal dilation coefficient (steel) & \si{K^{-1}} \\
$S_{\rm s}$ & Cable cross-sections & \si{m^2} \\
\end{tabular}
\end{center}
\caption{Parameters for the one-dimensional mechanical model (steel) \label{sec:THMmodel:subsec:consEq:subsubsec:THEq:table:parameters:1d}}
\end{figure}

Details and informations of the physical parameters for the three-dimensional are provided on Figure \ref{sec:THMmodel:subsec:consEq:subsubsec:THEq:table:parameters}, whereas those on the one-dimensional are given on Figure \ref{sec:THMmodel:subsec:consEq:subsubsec:THEq:table:parameters:1d}.

\subsection{Representative Structural Volume : standard section of a nuclear containment building\label{sec:THMmodel:subsec:RSV}}

The physical model is designed to capture the behavior of the so-called standard zone of the model, which corresponds to a portion of the mesh at mid-height, in the cylindrical part of the NCB. Thus, the region covered by the RSV comprises a three-dimensional portion containing three tangential prestressing cables and two vertical cables. For the section studied in this study, the internal radius of the wall is 21.9m, while the external radius is 23.4m. The width of the standard section corresponds to an angular sector of 4.2. For the scope of this work, the effect of passive steel reinforcement is neglected.

\begin{figure}[h!]
\begin{center}
\centering
\begin{subfigure}[b]{0.45\textwidth}
\begin{center}
\includegraphics[scale=0.35]{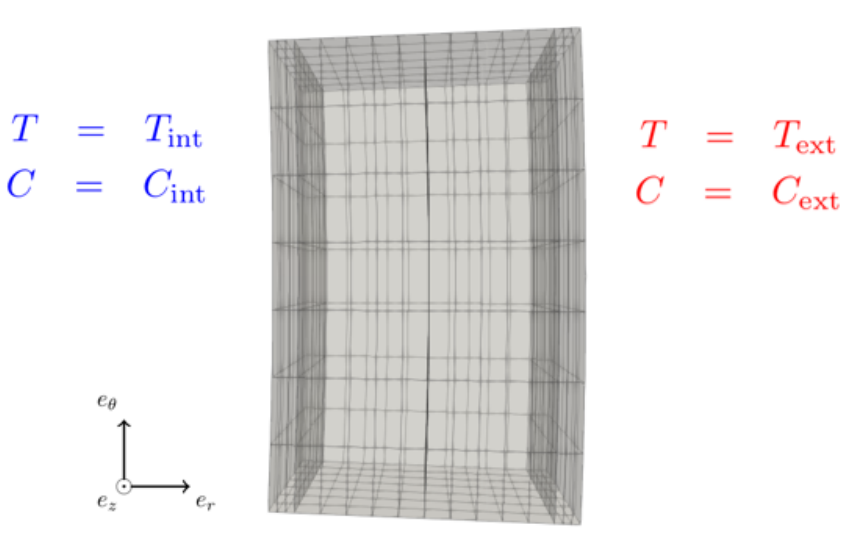} 
\end{center}
\caption{\normalsize  Temperature and water content BCs}
\end{subfigure}
\qquad
\centering
\begin{subfigure}[b]{0.45\textwidth}
\centering
\includegraphics[scale=0.65]{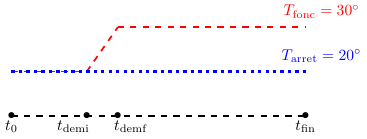} 
\includegraphics[scale=0.65]{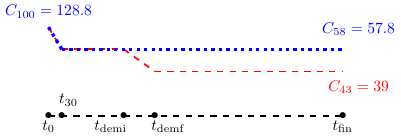}
\caption{\normalsize Temperature and water content evolutions}
\end{subfigure}
\end{center}
\caption{Boundary conditions for the thermal and hydric problems visualized on the HF thermal mesh \label{fig:thermalmehs_bcs}}
\end{figure}

Two mesh designs are used in practice: one for thermo-hydric calculations and another for mechanical calculations. The thermal mesh is refined close to the intrados and extrados to enable better reconstruction of the thermo-hydric gradients. The fields resulting from this procedure are then projected onto the mechanical mesh. The meshes employed in these studies are fairly coarse. In fact, these meshes have been built in order to be able to carry out uncertainty quantification or data assimilation studies. Therefore, engineers had to strike a balance between affordable computational cost and approximation quality.  Numerical solutions for thermal problems may exhibit oscillations (in terms of temporal and spatial discretizations). This may imply a violation of the maximum principle. To avoid this phenomenon, linear finite elements and a lumping of the mass matrix are used for this study. As previously mentioned, the thermal mesh does not contain the prestressing cables: it is composed of linear hexahedral cells (HEXA8). For the mechanical mesh, hexahedral quadratic elements (HEXA20) are employed for the concrete, and prestressing tendons are represented using SEG2 linear finite elements (2-node beams).

\begin{figure}[h!]
\begin{center}
\begin{subfigure}[b]{0.45\textwidth}
\begin{center}
\includegraphics[scale=0.2]{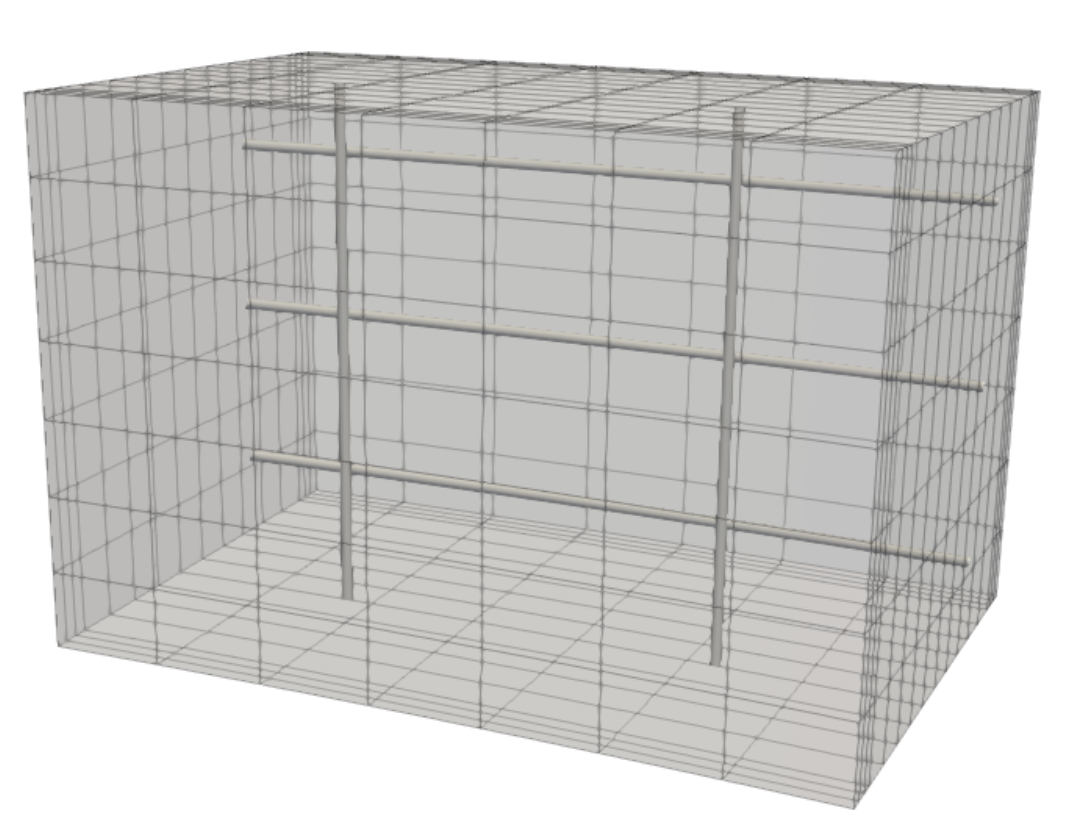} 
\end{center}
\caption{\normalsize Vizualisation of the mechanical mesh\label{fig:visu_infor_mecmesh:visu}}
\end{subfigure}
\qquad
\centering
\begin{subfigure}[b]{0.45\textwidth}
\begin{center}
\begin{adjustbox}{width=\columnwidth,center}
\begin{tabular}{ccccccc}
\hline
$N_{\rm e}$ & $N_{\rm e}^{1d}$ & $N_{\rm e}^{2d}$ & $N_{\rm e}^{3d}$ & $\mathcal{N}$ & $\mathcal{N}_{\rm c}$ & $\mathcal{N}_{\rm s}$ \\ \hline
1532 & 784 & 693 & 55 & 4076 & 3911 & 165
\end{tabular}
\end{adjustbox}
\caption{Summary of the parameters for the one-dimensional mechanical model\label{fig:visu_infor_mecmesh:info}}
\end{center}
\end{subfigure}
\end{center}
\caption{Visualization of the mechanical mesh (Figure \ref{fig:visu_infor_mecmesh:visu}) and information on the mechanical mesh (number of elements and number of nodes for one- and three-dimensional meshes, Figure \ref{fig:visu_infor_mecmesh:info})\label{fig:visu_infor_mecmesh}}
\end{figure}

The BCs and loads applied to the RSV zone are detailed below (Eq.\eqref{eq:thermalmehs_bcs}). Figure \ref{fig:thermalmehs_bcs} shows the temperature and water content histories adopted for the thermo-hydraulic calculations. As mentioned above, the BCs applied are Dirichlet conditions for temperature and water content. These are imposed on the inner wall (intrados) and the outer wall (extrados), as follows:

\begin{equation}
\label{eq:thermalmehs_bcs}
\left\{
\begin{array}{rclcl}
T &=& T_{\rm int}, &\text{on} & \Gamma_{\rm ext}, \\
T &=& T_{\rm ext}, &\text{on} & \Gamma_{\rm int}, \\
\end{array}
\right.
\quad \text{and} \quad
\left\{
\begin{array}{rclcl}
C &=& C_{\rm int}, &\text{on} & \Gamma_{\rm ext}, \\
C &=& C_{\rm ext}, &\text{on} & \Gamma_{\rm int}. \\
\end{array}
\right.
\end{equation}

\begin{figure}[h!]
\begin{center}
\includegraphics[scale=0.25]{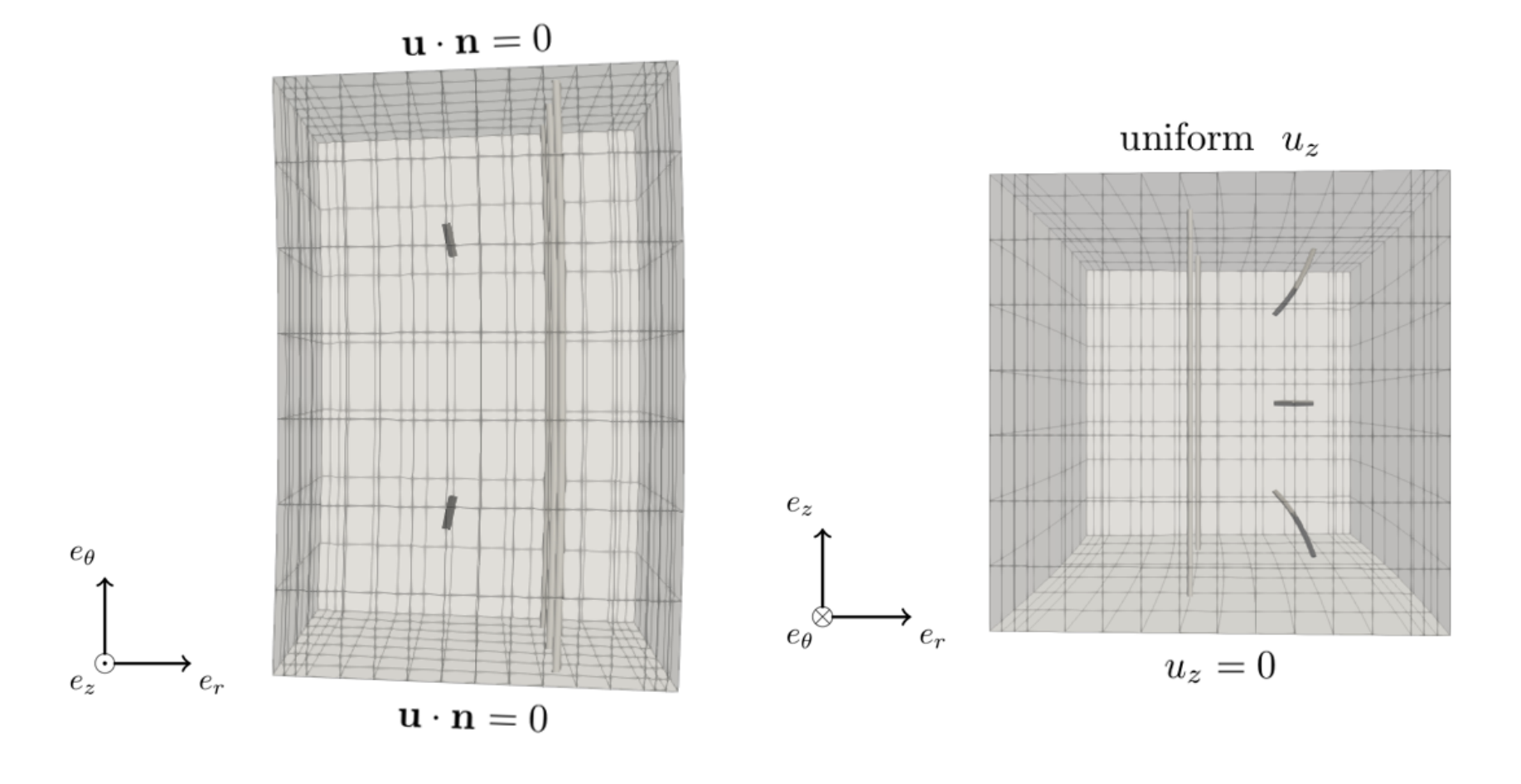} 
\end{center}
\caption{Boundary conditions for the mechanical problem visualized on the HF mechanical mesh\label{fig:mecmesh_bcs}}
\end{figure}

With regard to mechanical BCs, axisymmetric conditions are specified at the lateral boundaries of the RSV: this implies that normal displacements are assumed to be zero on each lateral face. Furthermore, vertical displacement is assumed to be blocked on the inner face of the RSV, while a uniform vertical displacement is used on the upper face. The set of boundary conditions with a visualization of the mechanical mesh is illustrated in Figure \ref{fig:mecmesh_bcs}.

\section{Numerical results: application to a standard section of a nuclear containment building}
\subsection{Solution Reproduction Problem\label{sec:numres_srpb}}

We first perform a validation of the methodology on a non-parametric case. We aim to mimic the HF simulation with our ROM for the same set of parameters. To assess the quality of the reduced model, we introduce several metrics. First of all, since our ROM is founded on a projection onto displacement modes, we introduce displacement approximation errors, at a given time step ($E^{{\rm app}, (k)}_{u, \mu}$), and averaged over time ($E^{{\rm app}, {\rm avg}}_{u, \mu}$):

\begin{equation}
\label{sec:numres_srpb:approxError}
E^{{\rm app}, (k)}_{u, \mu}= \frac{\norm{\fe{u}^{{\rm hf}, (k)}_\mu-\widehat{\fe{u}}^{(k)}_\mu}^2_2}{\norm{\fe{u}^{{\rm hf}, (k)}_\mu}^2_2}, \quad \text{and} \quad E^{{\rm app}, {\rm avg}}_{u, \mu}= \frac{\sqrt{\sum\limits_{k=1}^K\frac{t^{(k)} - t^{(k-1)}}{t_{\rm f}} \norm{\fe{u}^{{\rm hf}, (k)}_\mu-\widehat{\fe{u}}^{(k)}_\mu}^2_2}}{\sqrt{\sum\limits_{k=1}^K\frac{t^{(k)} - t^{(k-1)}}{t_{\rm f}} \norm{\fe{u}^{{\rm hf}, (k)}_\mu}^2_2}},
\end{equation}

\noindent where $t_{\rm f}$ is the final physical time used in the simulation and where $\fe{u}^{{\rm hf}, (k)}_\mu$ and $\widehat{\fe{u}}^{(k)}_\mu$ are respectively the solution at the k-th timestep obtained when using the HF model or the ROM for the parameter $\mu$. For the simulations reported below, we simulate a physical time of around 18 years.

\subsubsection{HF problem}

In this section, we present the HF problem we wish to reproduce. As previously stated, we are only seeking to reduce the mechanical calculation in our THM coupling. To this end, we rely on a thermo-hydraulic calculation, which can be viewed as an initial state common to all parametric calculations. These two simulations are carried out in compliance with the BCs described previously. On the figures provided afterwards, the time is given in seconds, as this is the time used in the numerical code (1 day = 86400 seconds). The time scheme for our creep simulations features an adaptive time step algorithm. In practice, in all the simulations carried out as part of this study, the entire simulation is performed over around 50 time steps.

\begin{figure}[h!]
\begin{center}
\begin{subfigure}[b]{0.45\textwidth}
\begin{center}
\includegraphics[scale=0.25]{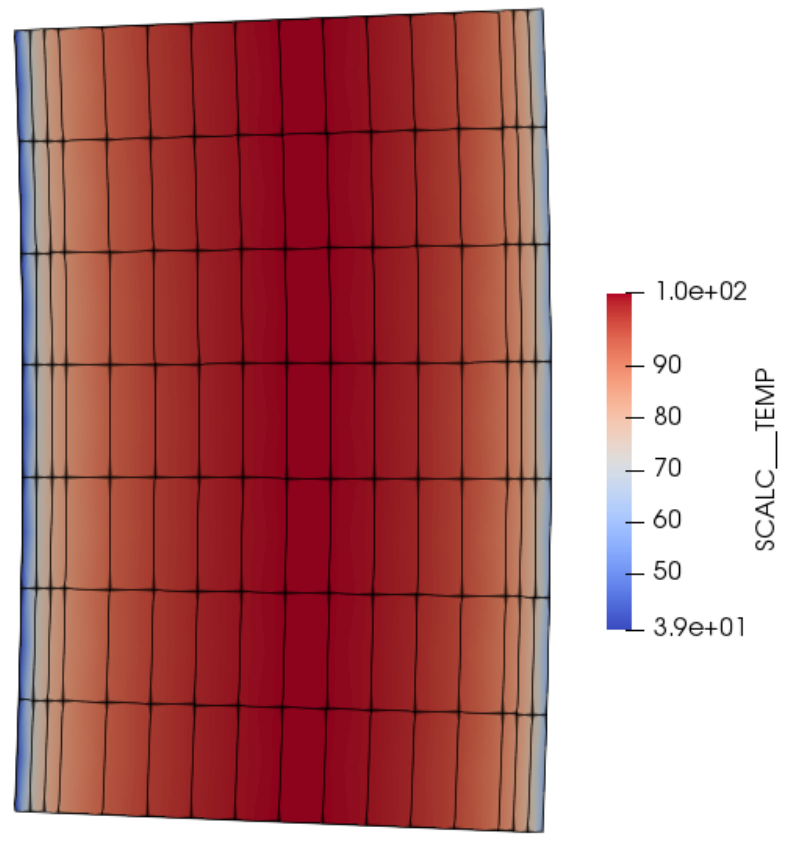}
\end{center}
\caption{\normalsize View of the drying field $C_w$ at the last time step of the HF simulation (top view)}
\end{subfigure}
\quad
\centering
\begin{subfigure}[b]{0.45\textwidth}
\begin{center}
\includegraphics[scale=0.775]{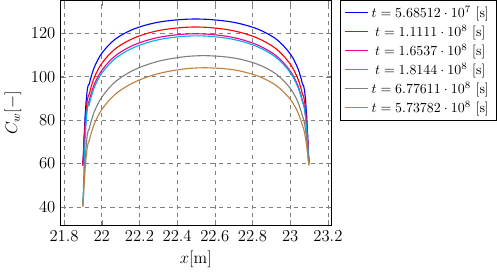} 
\end{center}
\caption{\normalsize Evolution of the drying field $C_w$ along $x$ in the plane ($y=0$, $z=0$)}
\end{subfigure}
\end{center}
\caption{Water content snapshots (output of the hydric calculation step) at the end of the HF simulation \label{sec:srpb_numres:subsec:hf:subsubsec:wc_field}}
\end{figure}

\begin{figure}[h!]
\begin{center}
\begin{subfigure}[b]{0.45\textwidth}
\begin{center}
\includegraphics[scale=0.25]{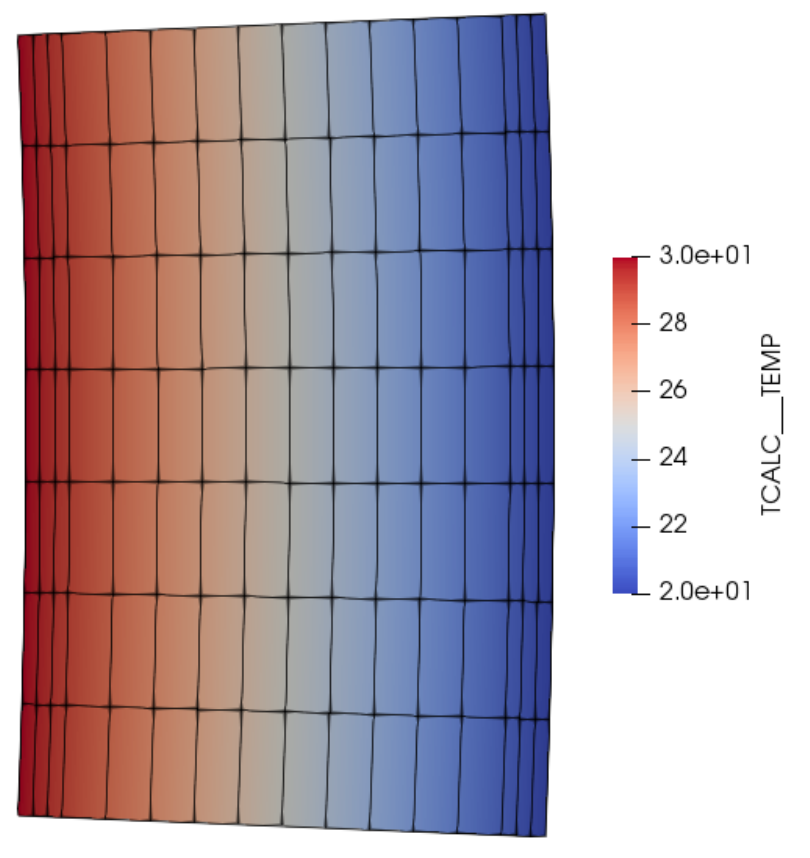}
\end{center}
\caption{\normalsize View of the temperature field $T$ at the last time step of the HF simulation (top view)}
\end{subfigure}
\quad
\centering
\begin{subfigure}[b]{0.45\textwidth}
\begin{center}
\includegraphics[scale=0.775]{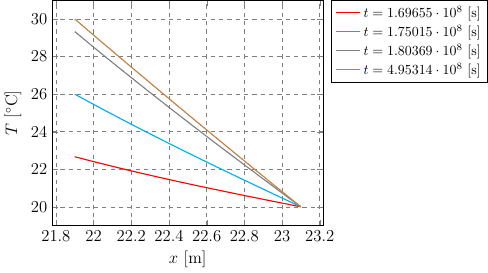} 
\end{center}
\caption{\normalsize  Evolution of the temperature field $T$ along $x$ in the plane ($y=0$, $z=0$)}
\end{subfigure}
\end{center}
\caption{Temperature snapshots  (output of the thermal calculation step) at the end of the HF simulation \label{sec:srpb_numres:subsec:hf:subsubsec:temp_field}}
\end{figure}

Figure \ref{sec:srpb_numres:subsec:hf:subsubsec:wc_field} displays the water content in the standard section at the end of the HF calculation. This figure depicts the evolution of the  $C_w$ field in the thickness of the containment building (in the standard section). Likewise, Figure \ref{sec:srpb_numres:subsec:hf:subsubsec:temp_field} shows the evolution of the temperature field in the thickness of the standard section. The physical parameters used for these calculations are summarized in Table \ref{sec:numres_param_table:coeff} where undefined parameters are chosen as follows:

$$
\overline{\eta}_{\rm dc}=5\cdot 10^{9}, \quad\overline{\kappa}=4.2\cdot 10^{-4}, \quad \overline{\alpha}_{\rm dc}=7.56\cdot 10^{-6}, \quad \overline{\eta}_{\rm is}=2.76\cdot 10^{18}, \quad \overline{\eta}_{\rm id}= 1.38\cdot 10^{18}.
$$ 
 
\begin{figure}[h!]
\begin{center}
\centering
\begin{subfigure}[b]{0.2\textwidth}
\begin{center}
\includegraphics[scale=0.215]{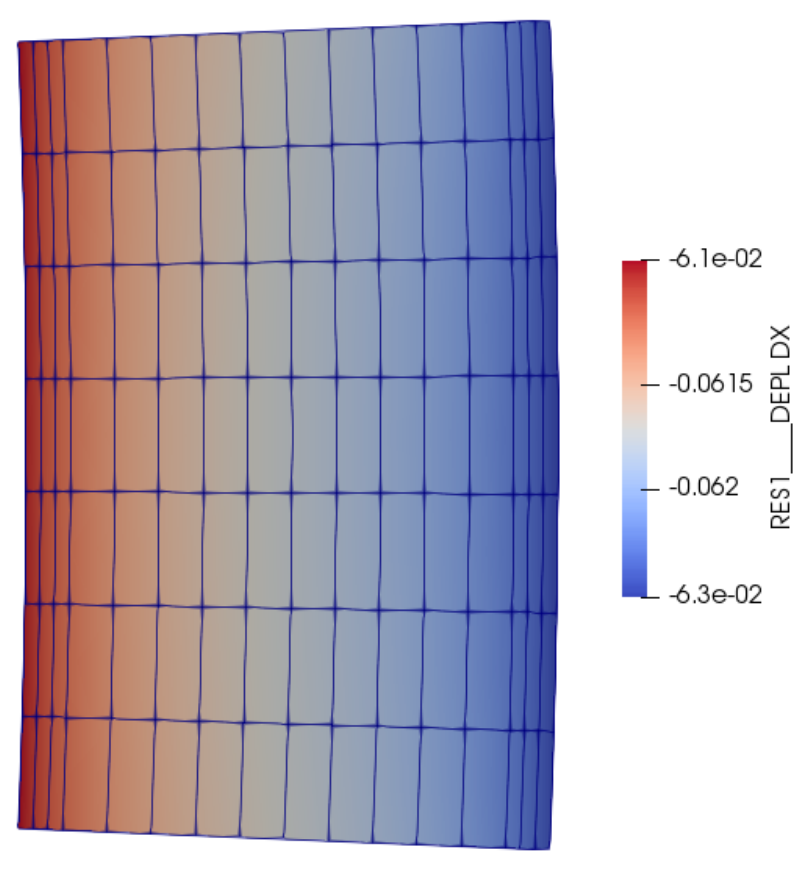} 
\end{center}
\caption{\normalsize $u_r$ $[\si{m}]$\label{sec:srpb_numres:subsec:hf:subsubsec:meca_field:ux}}
\end{subfigure}
\quad
\centering
\begin{subfigure}[b]{0.2\textwidth}
\begin{center}
\includegraphics[scale=0.215]{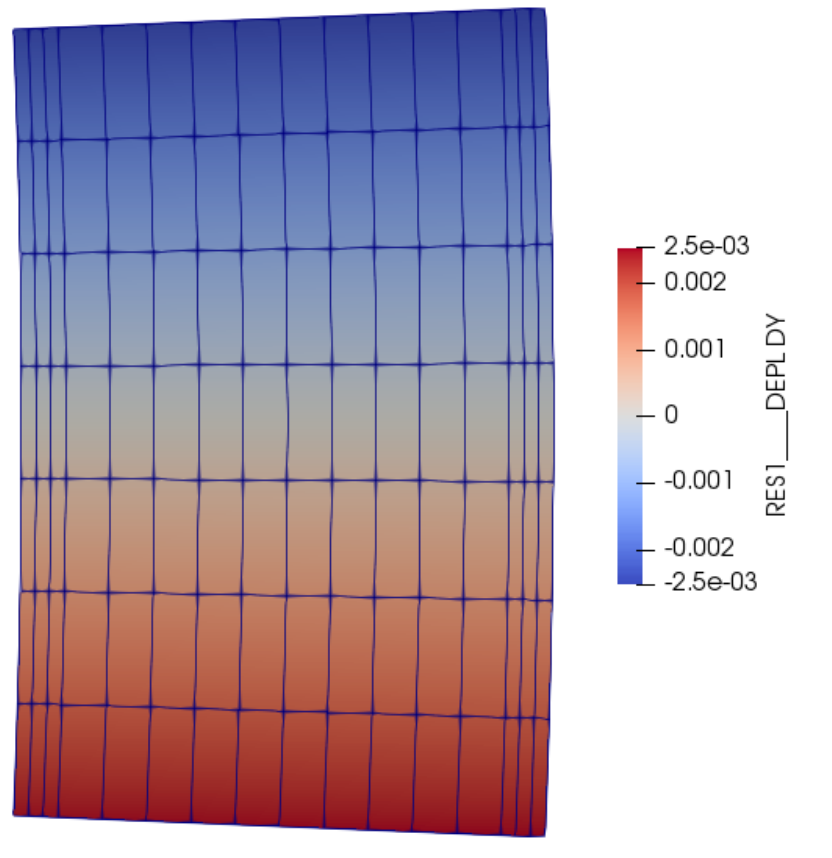} 
\end{center}
\caption{\normalsize$u_{\theta}$ $[\si{m}]$\label{sec:srpb_numres:subsec:hf:subsubsec:meca_field:uy}}
\end{subfigure}
\quad
\centering
\begin{subfigure}[b]{0.2\textwidth}
\begin{center}
\includegraphics[scale=0.175]{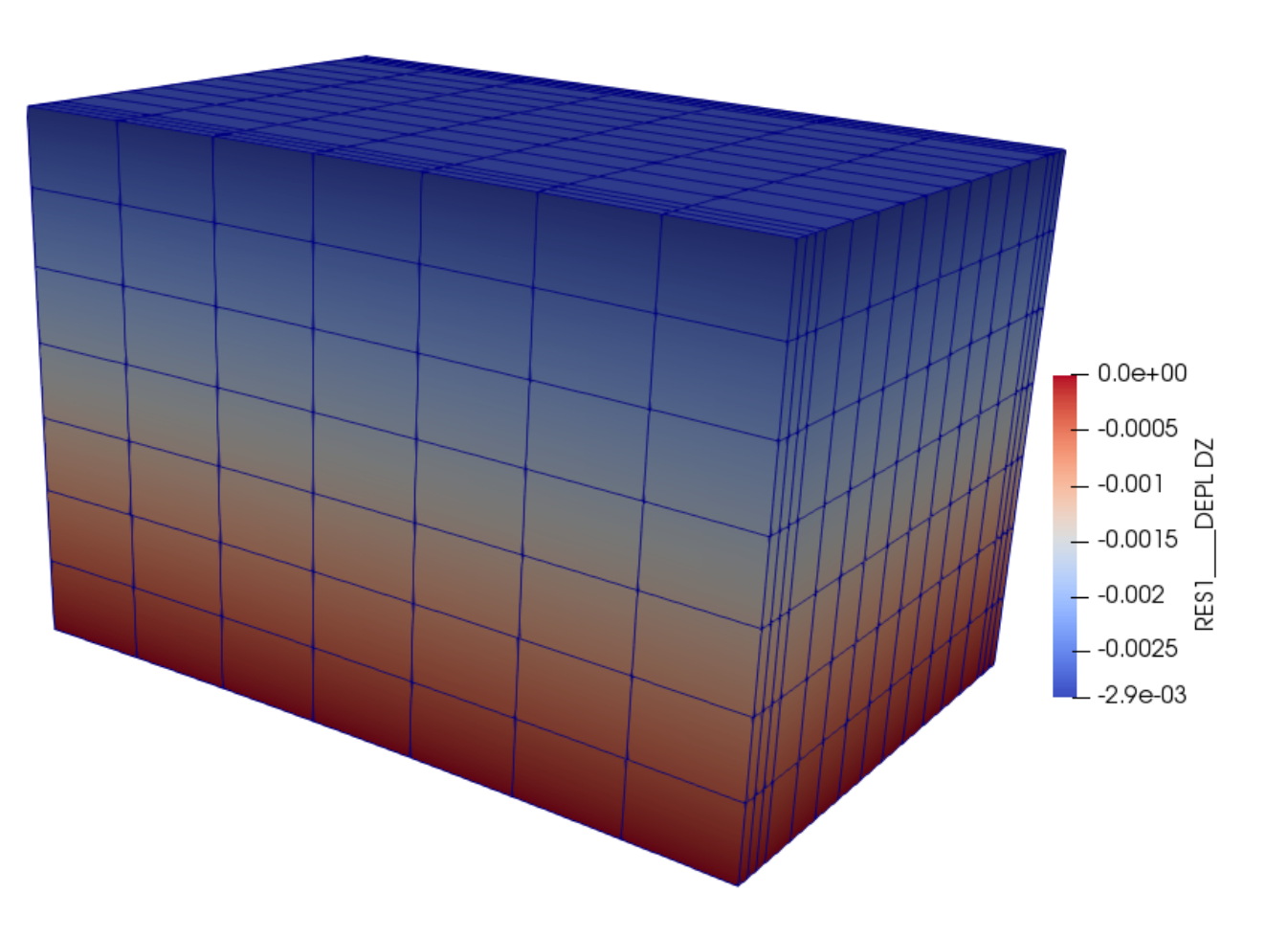} 
\end{center}
\caption{\normalsize $u_z$ $[\si{m}]$\label{sec:srpb_numres:subsec:hf:subsubsec:meca_field:uz}}
\end{subfigure}
\qquad
\centering
\begin{subfigure}[b]{0.2\textwidth}
\begin{center}
\includegraphics[scale=0.215]{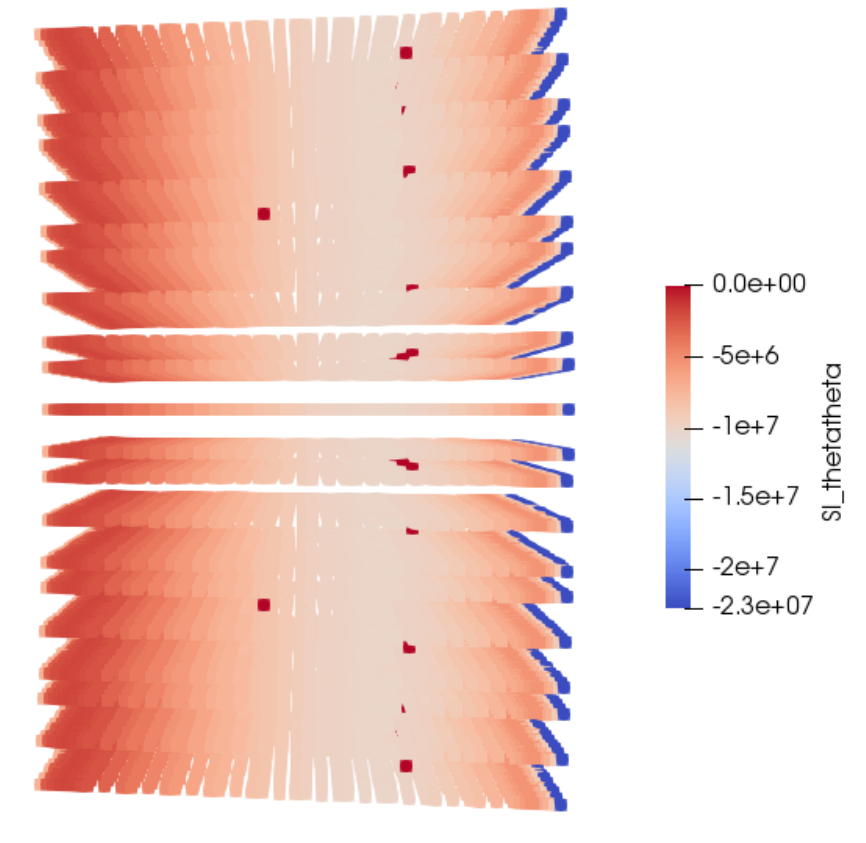} 
\end{center}
\caption{\normalsize \normalsize $\sigma_{\theta \theta}$ $[\si{Pa}]$\label{sec:srpb_numres:subsec:hf:subsubsec:meca_field:sigxx}}
\end{subfigure}
\end{center}
\caption{Mechanical fields snapshots (displacements, see Figure \ref{sec:srpb_numres:subsec:hf:subsubsec:meca_field:ux}, \ref{sec:srpb_numres:subsec:hf:subsubsec:meca_field:uy}, \ref{sec:srpb_numres:subsec:hf:subsubsec:meca_field:uz}, and stresses within the concrete, see Figure \ref{sec:srpb_numres:subsec:hf:subsubsec:meca_field:sigxx}) at the end of the HF calculation on the standard section\label{sec:srpb_numres:subsec:hf:subsubsec:meca_field}}
\end{figure}

 From these auxiliary fields ($H$ field in the methodology formulation in section \ref{sec:ROMmethodo}) we can determine all the mechanical fields using the HF code. Figure \ref{sec:srpb_numres:subsec:hf:subsubsec:meca_field} represents the displacement fields and the components of the Cauchy stress tensor obtained for the HF calculation we are seeking to reproduce in this section.
 
Our first goal is to ensure that the mechanical fields (displacements, stresses in the concrete and normal forces in the cables) are fairly accurate approximations of the values obtained from HF calculations. Besides, using a ROM of a standard section should provide a good quality approximation of the fields used in practical applications by engineers. In our case, this RSV has two main purposes: first, to compute leakage estimates from prestress loss in the cables, and second, to perform recalibration tests from deformation data (tangential and vertical deformations) on the intrados and extrados.

\begin{figure}[h!]
\begin{center}
\includegraphics[scale=0.85]{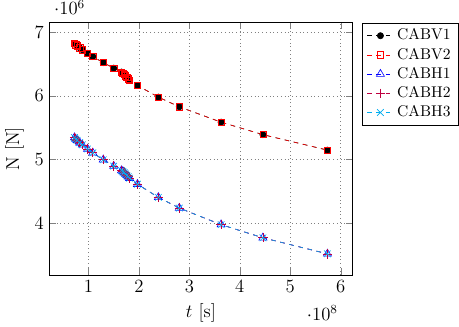} 
\end{center}
\caption{Evolution of normal forces in the two vertical (CABV1, CABV2) and three horizontal (CABH1, CABH2, CABH3) cables of the standard section\label{sec:srpb_numres:subsec:hf:subsubsec:meca_field:fig:Ncables}}
\end{figure}

\begin{figure}[h!]
\begin{center}
\centering
\begin{subfigure}[b]{0.45\textwidth}
\begin{center}
\includegraphics[scale=0.85]{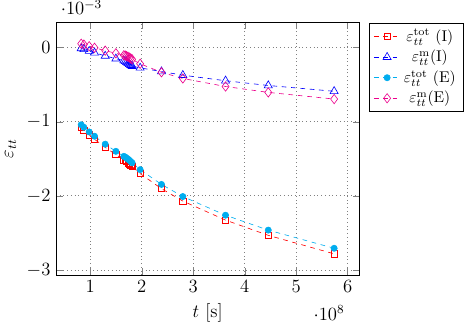} 
\end{center}
\caption{ \normalsize $\varepsilon_{tt}$}
\end{subfigure}
\qquad
\centering
\begin{subfigure}[b]{0.45\textwidth}
\begin{center}
\includegraphics[scale=0.85]{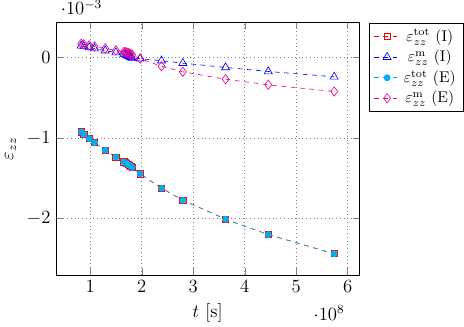} 
\end{center}
\caption{ \normalsize $\varepsilon_{zz}$}
\end{subfigure}
\end{center}
\caption{Comparison for the pointwise values between some components (tangential and vertical) the mechanical strains and the total strains in sensor zones (extrados (E) and intrados (I)) \label{sec:srpb_numres:subsec:hf:subsubsec:meca_field:fig:defpointwise}}
\end{figure}

Figure \ref{sec:srpb_numres:subsec:hf:subsubsec:meca_field:fig:Ncables} depicts the evolution of the mean value of the normal forces in each of the five cables within the standard section. For the record, the mesh studied contains two vertical cables and three horizontal cables. Within the framework of the investigated model, the vertical cables have a similar behavior (as do the three horizontal cables). In the following, we have decided to report only the results for one horizontal and one vertical cables (CABV1 and CABH2), to ease the readability of the results. Figure \ref{sec:srpb_numres:subsec:hf:subsubsec:meca_field:fig:defpointwise} displays the evolution of mechanical strains and total strains in the concrete. In our notations, (I) stands for intrados whereas (E) stands for extrados. In our cases of interest, the total strains of the material are not purely mechanical. In general, data assimilation problems only focus on mechanical deformations. This is of key interest when reconstructing the strain field from the displacement modes, since the strain includes components due to temperature gradients and/or water pressure. Indeed, in our ROM resolution procedure, we have generalized coordinates at our disposal, which enable us to reconstruct the displacement field in the material. By computing the symmetric gradient of this displacement field, we can determine the total strains. In order to reconstruct a strain field, we must subtract the terms related to the thermal and hydric fields. Both these fields may be derived independently of the reduction process, since we only reduce the mechanical part of the calculation chain. We are thus able to pre-calculate the TH strain fields and subtract them from a total strain field so as to obtain the reconstructed mechanical strain field.

In order to assess the accuracy of our reduced model, we introduce approximation errors for these different fields: for the average of the normal forces at the nodes in the CABV1 vertical cable ($E_{\overline{\mu}}^{{\text{app}}, (t)} [\text{N}_{\rm V_2}]$), and in the horizontal cable ($E_{\overline{\mu}}^{{\text{app}}, (t)} [\text{N}_{\rm H_2}]$), for the average of the tangential strain and vertical strain on the extrados ($E_{\overline{\mu}}^{{\text{app}}, (t)}[\varepsilon^{\rm m}_{tt} \text{ (avg - E) }]$ and $E_{\overline{\mu}}^{{\text{app}}, (t)}[\varepsilon^{\rm m}_{zz} \text{ (avg - E) }]$) , and finally for the average of the tangential strain and horizontal strain on the intrados ($E_{\overline{\mu}}^{{\text{app}}, (t)}[\varepsilon^{\rm m}_{tt} \text{ (avg - I) }]$ and $E_{\overline{\mu}}^{{\text{app}}, (t)}[\varepsilon^{\rm m}_{zz} \text{ (avg - I) }]$) . To average the components of the strain tensor, the values at the Gauss points are extrapolated to the nodes, and the value at the nodes is then averaged. These relative errors in the deformation fields relate exclusively to mechanical deformations. Indeed, this is the only part of the tensor that is actually modified by our reduction process, as explained above.

\subsubsection{Speedups and approximation errors}

In order to validate the ROM, we verify that the displacement field is properly reconstructed. Furthermore, since we are interested in the use of the ROM for engineering applications, it is necessary to confirm the quality of the approximation on the various quantities of interest, more precisely tangential and vertical deformations and normal forces in the cables (which enables us to calculate prestressing loss). Ultimately, it is crucial to provide a model that reduces the computation time required whenever a call is made. To this end, we focus on the speedups ($\text{speedup}=\frac{\text{ROM CPU cost}}{\text{HF CPU cost}}$) obtained after construction of the reduced model (online phase).

\begin{figure}[h!]
\begin{center}
\includegraphics[scale=0.85]{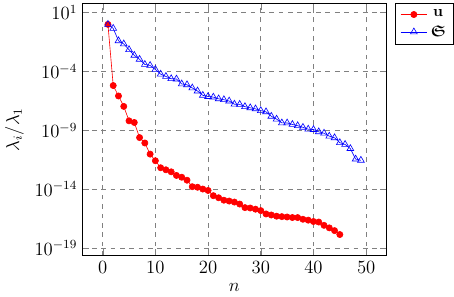} 
\caption{POD eigenvalues for the displacement ($\fe{u}$) and the generalized forces ($\bm{\mathfrak{S}}$) using a $\ell_2$ compression\label{sec:numres:subsec:srpb:fig:eigenvales_ell2} for a solution reproduction problem (50 initial snapshots)}
\end{center}
\end{figure}

Figure \ref{sec:numres:subsec:srpb:fig:eigenvales_ell2} depicts the POD eigenvalues generated on snapshots of displacements ($\fe{u}$) and generalized forces ($\bm{\mathfrak{S}}$). The decay profiles are quite distinct between the two physical quantities: the decay of the eigenvalues for displacements is fast, unlike in the case of generalized forces. This implies that the sizes of the two bases generated for POD tolerances of the same order of magnitude are significantly different. The displacement basis will always be much smaller than the generalized force basis.

\begin{figure}[h!]
\begin{center}
\centering
\begin{subfigure}[b]{0.45\textwidth}
\begin{center}
\includegraphics[scale=0.85]{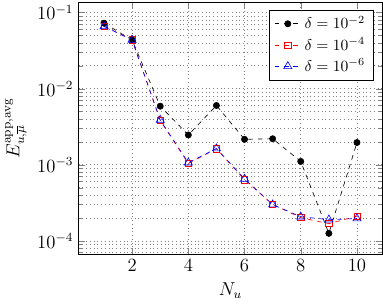} 
\caption{Approximation errors on $\fe{u}$\label{numres:srpb:approxerrorsu_speedups:ae}}
\end{center}
\end{subfigure}
\quad
\centering
\begin{subfigure}[b]{0.45\textwidth}
\begin{center}
\includegraphics[scale=0.85]{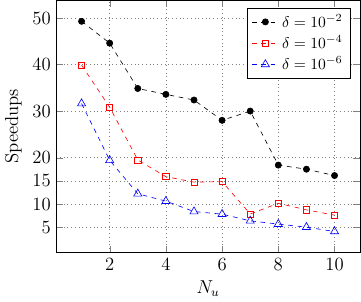} 
\caption{Speedups\label{numres:srpb:approxerrorsu_speedups:sp}}
\end{center}
\end{subfigure}
\end{center}
\caption{Evolution of time-averaged approximation errors on the displacements and speedups as a function of the number of modes used ($N_u$, see Figure\ref{numres:srpb:approxerrorsu_speedups:ae}) and for several hyper-reduction tolerances ($\delta$, see Figure\ref{numres:srpb:approxerrorsu_speedups:sp})\label{numres:srpb:approxerrorsu_speedups}}
\end{figure}

\begin{figure}[h!]
\begin{center}
\centering
\begin{subfigure}[b]{0.29\textwidth}
\begin{center}
\includegraphics[scale=0.2]{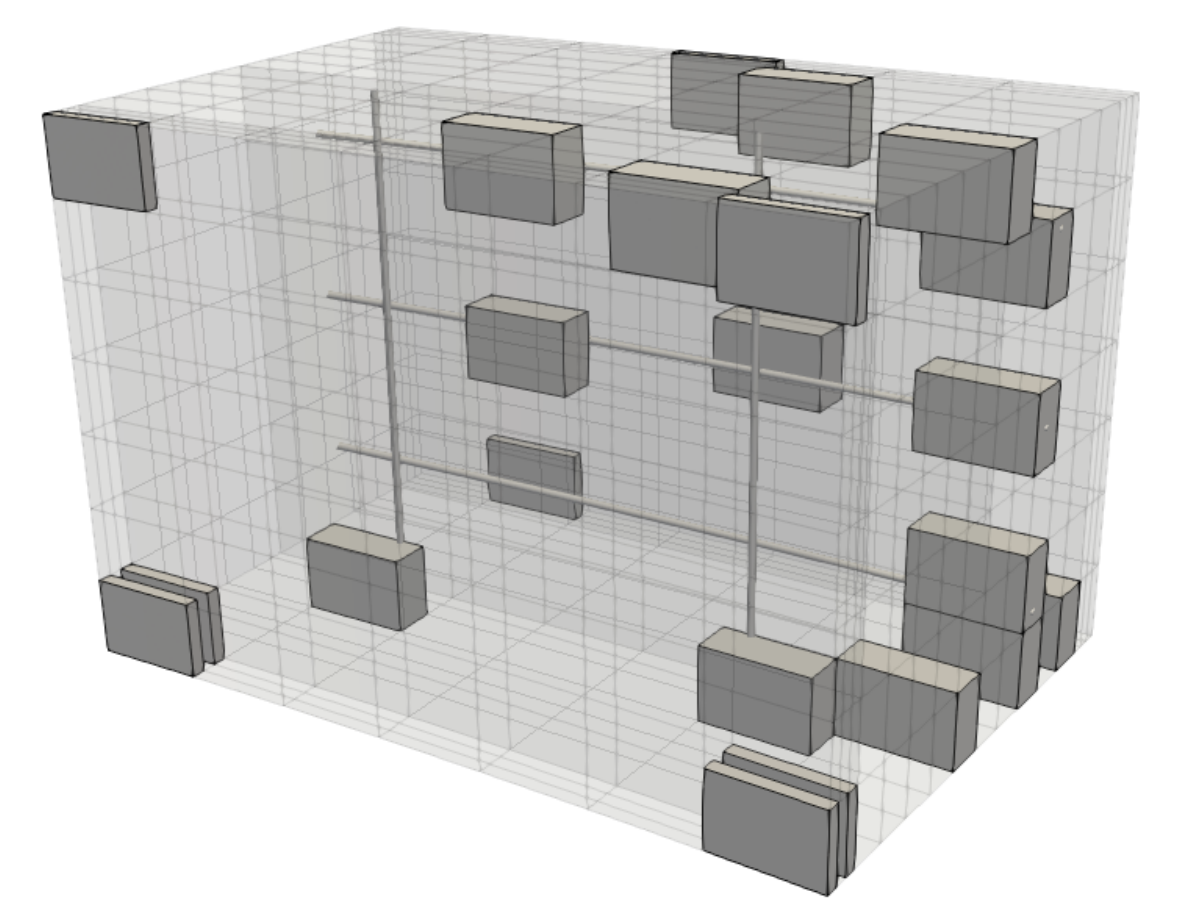}
\end{center}
\caption{\normalsize $\delta=10^{-2}$}
\end{subfigure}
\quad
\centering
\begin{subfigure}[b]{0.29\textwidth}
\begin{center}
\includegraphics[scale=0.2]{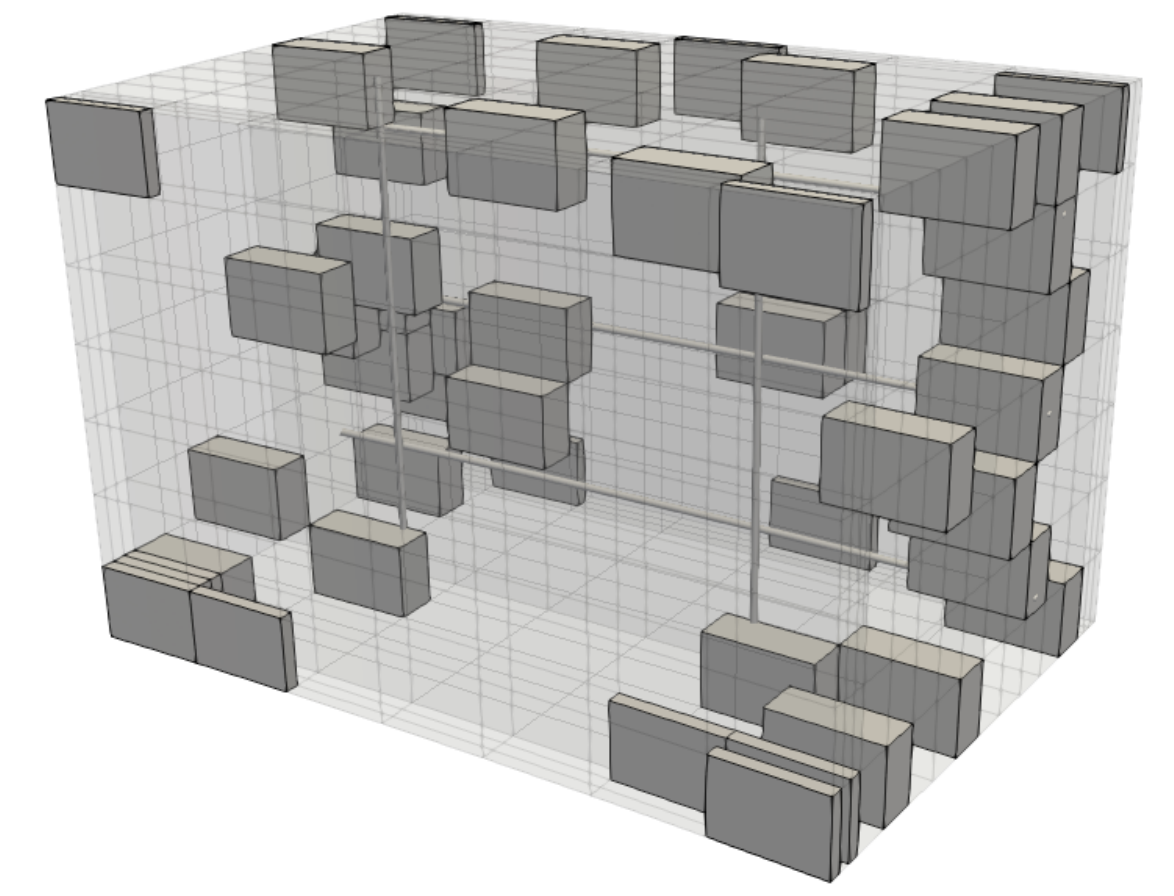}
\end{center}
\caption{\normalsize $\delta=10^{-4}$}
\end{subfigure}
\quad
\centering
\begin{subfigure}[b]{0.29\textwidth}
\begin{center}
\includegraphics[scale=0.2]{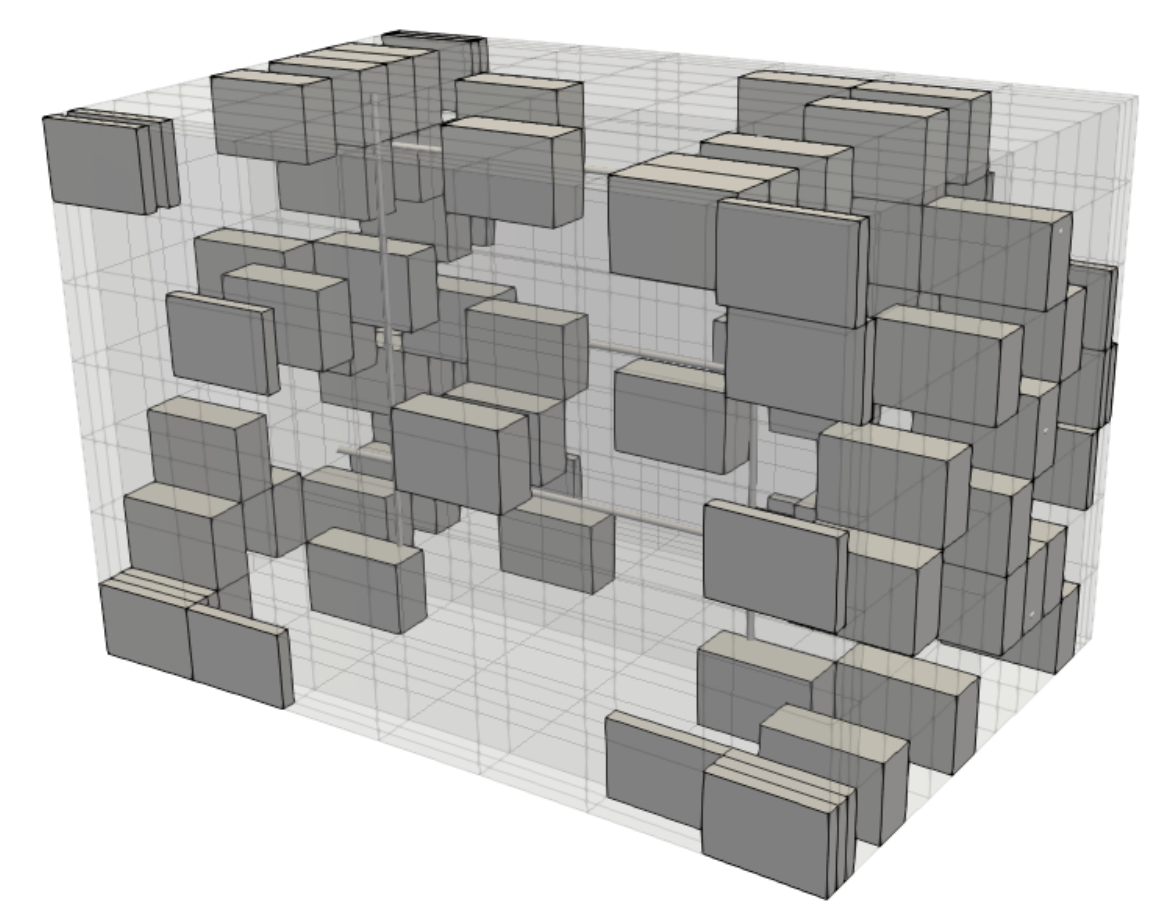}
\end{center}
\caption{\normalsize $\delta=10^{-6}$}
\end{subfigure}
\end{center}
\caption{Reduced meshes of the standard section obtained for a reproduction problem solution using $N_u=5$ displacement modes and for several hyper-reduction parameters \label{numres:srpb:reducedmeshes}}
\end{figure}

As a way of assessing the robustness of the reduction approach proposed here, we have built several ROMs for different numbers of displacement modes and different hyper-reduction tolerances. An increase in the number of modes and a decrease in the $\delta$ hyperparameter both improve the quality of the ROM and increase computation time (speedup). Thus, a tradeoff needs to be found for engineering applications in order to provide a fast and accurate ROM. Figure \ref{numres:srpb:approxerrorsu_speedups} displays the evolution of speedups and time-averaged displacement approximation errors as a function of the number of modes (for several tolerances). We observe that from 5 modes upwards, The reduced order model exhibits an good approximation quality, with approximation errors below the order of 0.2$\%$ (for all tolerances studied). In this case, the speedups achieved are substantial: around 10 for the most severe tolerance (equal to the Newton-Raphson tolerance), around 15 for the intermediate tolerance studied, and over 30 for the coarsest tolerance. These accelerations in CPU computation time are all the more appealing as the mesh studied in this paper is very coarse, with only a few hundred elements (see Figure \ref{numres:srpb:reducedmeshes} for further details). This opens the door to future work on the use of finer meshes in NCB cross-section studies.

\begin{figure}[h!]
\begin{center}
\begin{subfigure}[b]{0.45\textwidth}
\begin{center}
\includegraphics[scale=0.85]{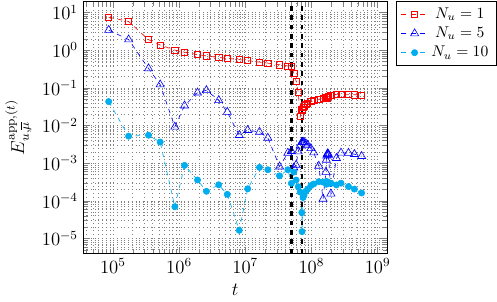} 
\end{center}
\caption{Time evolution of relative errors for $\delta=10^{-6}$ vectors in the reduced basis and varying number of $N_u$ values}
\end{subfigure}
\quad
\begin{subfigure}[b]{0.45\textwidth}
\begin{center}
\includegraphics[scale=0.85]{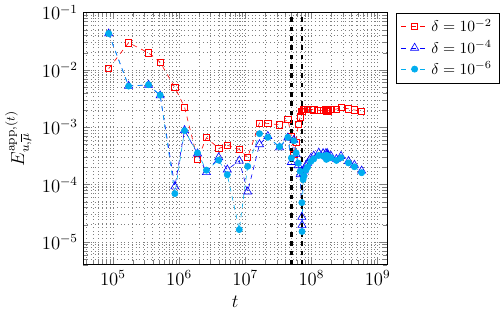} 
\end{center}
\caption{Time evolution of relative errors for $N_u=10$ vectors in the reduced basis and varying number of $\delta$ values}
\end{subfigure}
\end{center}
\caption{Evolution of approximation errors on displacements at each time step for several numbers of modes used or for several hyper-reduction tolerances\label{numres:srpb:approxerror_timeevol}}
\end{figure}

We have further investigated the quality of the ROM along the time trajectory of the problem. Figure \ref{numres:srpb:approxerror_timeevol} represents the relative errors at each time step for different ROMs. Since the construction of the ROM is determined by a pair of hyperparameters $\pr{N_u, \delta}$, we focused on the influence of each parameter in fixing the second. The parameters set in the two test cases are chosen so as to be as restrictive as possible in the parameter sets we explore here. We find that for our problem, the number of modes has a much greater influence on time-evolution profiles than hyper-reduction tolerance. Since the latter parameter leads to an increase in mesh size as it decreases, this prompts us to state that: in this non-parametric case, it is advisable to fix a number of modes to control the approximation error, and it suffices to take a low or intermediate tolerance to get good speedups. We notice that for low approximation qualities, there are jumps in the relative error profiles of the displacement fields. This is due to the fact that the ROM is built over the entire life of the standard section, namely with three distinct physical regimes: life of concrete without cables, prestressing, and life of concrete with cables. For small numbers of modes, the ROMs is unable to generate modes designed to approximate these three phases. Since we chose to use no weighting, it will have a tendency to approximate the final step much more accurately, which is justified by the fact that the number of time steps associated with this phase is much greater. This higher approximation quality on the last step is of interest for our applications, as we seek not only a reliable approximation in terms of time trajectories, but also, and above all, a solution that is truly representative of the system's final state. If we need control the time-averaged approximation errors in a different manner, it would be natural to use a weighted POD in order to take into account the non-constant timestepping.

\subsubsection{Errors on the quantities of interest}
\begin{figure}[h!]
\begin{center}
\begin{subfigure}[b]{0.31\textwidth}
\begin{center}
\includegraphics[scale=1]{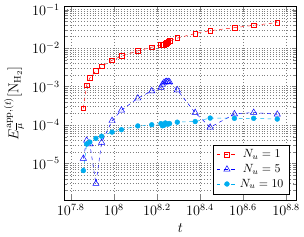} 
\end{center}
\end{subfigure}
\quad
\begin{subfigure}[b]{0.31\textwidth}
\begin{center}
\includegraphics[scale=1]{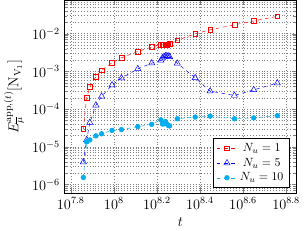} 
\end{center}
\end{subfigure}
\quad
\begin{subfigure}[b]{0.31\textwidth}
\begin{center}
\includegraphics[scale=1]{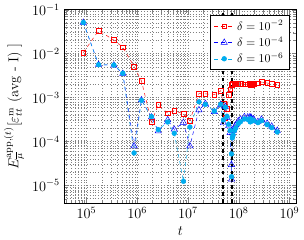} 
\end{center}
\end{subfigure}\\
\begin{subfigure}[b]{0.31\textwidth}
\begin{center}
\includegraphics[scale=1]{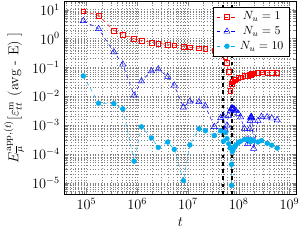} 
\end{center}
\end{subfigure}
\quad
\begin{subfigure}[b]{0.31\textwidth}
\begin{center}
\includegraphics[scale=1]{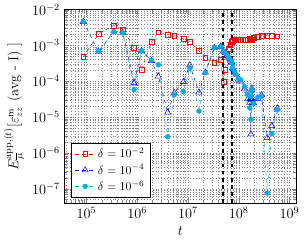} 
\end{center}
\end{subfigure}
\quad
\begin{subfigure}[b]{0.31\textwidth}
\begin{center}
\includegraphics[scale=1]{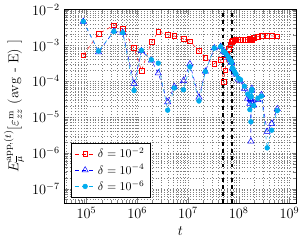}
\end{center}
\end{subfigure}
\end{center}
\caption{Evolution of approximation errors on QoIs at each time step for several numbers of modes used or for several hyper-reduction tolerances (the two vertical lines in black delimit the prestressing section of the cables)\label{numres:srpb:approxerror_qoIs_timeevol}}
\end{figure}

The scope of the research we have undertaken requires us to be confident in our ability to provide accurate QoIs. We thus wish to verify that the ROM obtained, in addition to being a good approximation of the HF calculation in terms of displacements while being significantly less computationally expensive, can be used in real applications. This is achieved by investigating the profiles of normal forces in the cables and deformations at the sensor level (average measure of a component of the strains tensor over the internal or external face). We would like to point out that data post-processing differs according to the QoIs studied. The reduced mesh contains all the prestressing cables, while the quadrature laws are unchanged in the one-dimensional mesh. As a result, we can compute the relative error on normal forces directly after calling up the reduced model. For strains, however, we must reconstruct the strain fields on the HF mesh, and then apply the observation operators (physical sensors) used in the HF framework. This step s computationally inexpensive compared to the overall procedure, as the symmetric gradients of the modes are already known, because they are required for the hyper-reduction process. All that needs to be done is to multiply these modes to the generalized coordinates and apply the observation operator. Figure \ref{numres:srpb:approxerror_qoIs_timeevol} provides the time-evolution of the relative errors on the QoIs. On Figure \ref{numres:srpb:approxerror_qoIs_timeevol}, we delimit the three phases of a mechanical calculation for a power plant containment building: a first phase in which the cables are not involved in the mechanical calculation, i.e. the concrete evolves on its own; a second phase in which the concrete is prestressed (see Eq. \eqref{eq:loadings_for_prestressed_step:anal} for specific loads in this case); then, finally, the life of the prestressed concrete, in which the concrete and cables are fully coupled. The three periods are delimited by dotted black vertical lines. The HF solver's adaptive time-stepping process explains the temporal distribution of the various snapshots. The initial time for plotting corresponds to the first time step output by the reference calculation code.

The pattern of strain changes is similar to that of displacement approximation errors. Furthermore, the observation of a better approximation of deformations during the life of the NCB after prestressing is also confirmed. This confirms the usefulness of the ROM for data assimilation problems. In practice, data is only available once the cables have been prestressed. For the sake of clarity, we would like to point out that the time scale for the profile of relative errors in normal forces is not the same as that for deformations. In fact, only the life of the enclosure after prestressing is depicted, since normal forces are always zero beforehand, or known analytically.

\subsection{Parametric problem\label{sec:numres_param}}

In a second step, we study a parametric case. As mentioned above, we consider here a strong-greedy approach. Thus, in order to drive the greedy search, we consider the maximum approximation error on a given training set ($\Theta_{\rm train}$), for the parameters we have not yet examined. As a reminder, $\Theta_*$ corresponds to the set of parameters used in building the ROM. We introduce a notation for the maximal error obtained when testing the ROM:

\begin{equation*}
\Delta_N^{\rm stg} = \max\limits_{i\in \Theta_{\rm train}\setminus\Theta_{*}}E^{{\rm app}, {\rm avg}}_{u, \mu_i}.
\end{equation*}

In the physical case under study, uncertainty is mainly limited to five physical parameters $\mu=[\eta_{\rm dc}, \kappa,, \alpha_{\rm dc}, \eta_{\rm is}, \eta_{\rm id}]^\top\in \R^5$, and in particular to the first two. As a validation of our model reduction approach, we set all the other parameters of the problem (see values in the Table \ref{sec:numres_param_table:coeff}), and restrict the parametric problem to the other parameters.

\begin{table}[h!]
\begin{center}
\begin{tabular}{ccccc}
\hline
\textbf{Input parameter} & \textbf{Notation} & \textbf{Value} & \textbf{Unit} \\ \hline
 Young's modulus (steel) & $E_{\rm s}$ & $1.9\cdot 10^{11}$ & \si{Pa} \\
 Poisson's ratio (steel) & $\nu_{\rm s}$ & $0.3$ & \si{-} \\
 Density (steel) & $\rho_{\rm s}$ & $7850$ & \si{kg.m^3} \\
 Thermal dilation coefficient (steel) & $\alpha_{\rm th, s}$ & $1\cdot 10^{-5}$ &\si{K^{-1}} \\
 Guaranteed maximum load stress at break & $f_{\rm prg}$ & $1.86\cdot 10^9$ & \si{Pa}\\
 Cable cross-section & $S_{\rm s}$ & $5400\cdot 10^{-6}$ & \si{m}\\
\hdashline
 Young's modulus (concrete) & $E_{\rm c}$ & $4.2\cdot 10^{10}$ & \si{Pa} \\
 Poisson's ratio (concrete) & $\nu_{\rm c}$ & $0.2$ & \si{-} \\
 Density (concrete) & $\rho_{\rm c}$ & $2350$ & \si{kg.m^3} \\
 Thermal dilation coefficient (concrete) & $\alpha_{\rm th, c}$ & $5.2\cdot 10^{-6}$ &\si{K^{-1}} \\
 Autogenous  shrinkage coefficient & $\beta_{\rm endo}$ & $66.1\cdot 10^{-6}$& \si{-} \\
 Dessication shrinkage coefficient & $\alpha_{\rm dc}$ & \textbf{X} & \si{-} \\
 Reversible deviatoric basic stiffness & $k_{\rm rd}$ &$5.98\cdot 10^{18}$ & \si{Pa} \\
 Reversible deviatoric basic viscosity & $\eta_{\rm rd}$ & $8.12\cdot 10^{16}$ &\si{Pa.s} \\
 Irreversible deviatoric basic viscosity & $\eta_{\rm id}$  & \textbf{X} &\si{Pa.s} \\
 Basic creep activation energy & $U_{\rm bc} / R$ & $4700$ &\si{K} \\
 Basic creep reference temperature & $T_{\rm bc}^0$ & $20$ &\si{^\circ C} \\
 Basic creep consolidation parameter & $\kappa$ & \textbf{X} &\si{-} \\
 Desiccation creep viscosity & $\eta_{\rm dc}$ & \textbf{X} &\si{Pa^{-1}} \\ \hdashline
 Dead weight of upper concrete lifts  & $\sigma_{z,c}$ & $1.375\cdot 10^6$ &\si{Pa} \\ \hdashline
 Stress applied to vertical cables & $\sigma_{v, s}$  & $990.7\cdot 10^6$ & \si{Pa} \\
 Stress applied to horizontal cables & $\sigma_{h, s}$ & $1264.7\cdot 10^6$ & \si{Pa} \\ \hline
\end{tabular}
\end{center}
\caption{Coefficients for the mechanical model fixed for the parametric problem. The notation \textbf{X} corresponds to the parameters that can vary and, therefore, we do not give a priori numerical values.\label{sec:numres_param_table:coeff}}
\end{table}

\subsubsection{In-sample test for $\mathcal{P}\subset\R^2$}
We confine the study to a parametric case with two parameters. The vector of parameters considered is as follows:

\begin{equation*}
\mu = \begin{bmatrix}
\eta_{\rm dc} \\
\kappa
\end{bmatrix}\in \left[5\cdot 10^8, \ 5 \cdot 10^{10}\right]\times\left[10^{-5}, \ 10^{-3}\right] \subset \R^2.
\end{equation*}

\noindent This is tantamount to setting the following parameters (in addition to those given in the Table \ref{sec:numres_param_table:coeff}):

$$
\overline{\alpha}_{\rm dc}=7.56\cdot 10^{-6}, \quad \overline{\eta}_{\rm is}=2.76\cdot 10^{18}, \quad \overline{\eta}_{\rm id}= 1.38\cdot 10^{18}.
$$

\begin{figure}[h!]
\begin{center}
\centering
\begin{subfigure}[b]{0.45\textwidth}
\begin{center}
\includegraphics[scale=0.85]{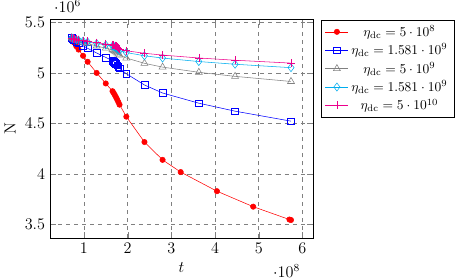}
\end{center}
\caption{ \normalsize $\kappa=1\cdot 10^{-5}$ (horizontal cables) \label{sec:numres_param:subsec:parametric2:fig:normal_forces:t11}}
\end{subfigure}
\qquad
\centering
\begin{subfigure}[b]{0.45\textwidth}
\begin{center}
\includegraphics[scale=0.85]{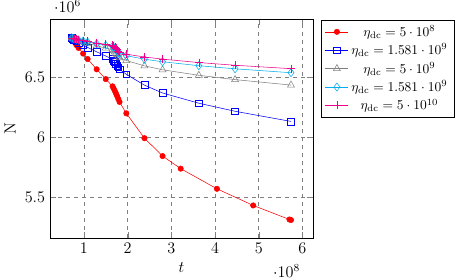}
\end{center}
\caption{ \normalsize $\kappa=1\cdot 10^{-5}$ (vertical cables)\label{sec:numres_param:subsec:parametric2:fig:normal_forces:t12}}
\end{subfigure}
\\
\centering
\begin{subfigure}[b]{0.45\textwidth}
\begin{center}
\includegraphics[scale=0.85]{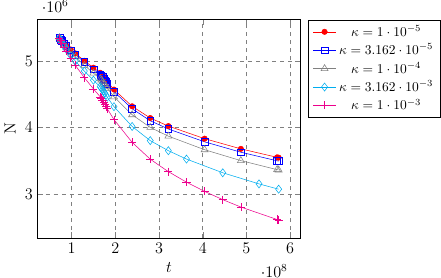}
\end{center} 
\caption{ \normalsize $\eta_{\rm dc} = 5\cdot 10^7$ (horizontal cables)\label{sec:numres_param:subsec:parametric2:fig:normal_forces:t21}}
\end{subfigure}
\qquad
\centering
\begin{subfigure}[b]{0.45\textwidth}
\begin{center}
\includegraphics[scale=0.85]{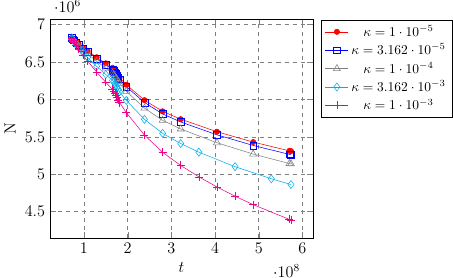}
\end{center}
\caption{ \normalsize $\eta_{\rm dc} = 5\cdot 10^7$ (vertical cables) \label{sec:numres_param:subsec:parametric2:fig:normal_forces:t22}}
\end{subfigure}
\end{center}
\caption{Evolution of normal forces over time for pairs of parameters belonging to the parametric set of size $|\Theta_{\rm train}|=25$. Figures \ref{sec:numres_param:subsec:parametric2:fig:normal_forces:t11}-\ref{sec:numres_param:subsec:parametric2:fig:normal_forces:t12} (resp. Figure \ref{sec:numres_param:subsec:parametric2:fig:normal_forces:t21}-\ref{sec:numres_param:subsec:parametric2:fig:normal_forces:t22}) feature cases where the parameter $\kappa$ (resp. $\eta_{\rm dc}$) is fixed. For each pair, we plot the time evolution of the normal forces averaged over all the nodes of the vertical and horizontal cables.\label{sec:numres_param:subsec:parametric2:fig:normal_forces}}
\end{figure}

We rely on a training space of size $|\Theta_{\rm train}|=25$, designed as the tensor product of two one-dimensional grids log-evenly spaced ($5\times 5$ grid). This choice results from a tradeoff between the need for sufficiently fine discretization to have several parameters, and the offline CPU cost of building the ROM (an HF calculation takes around fifteen minutes). The choice of optimal discretization is out of the scope of this work and is a field of research of its own. To help understand the physical problem under study, Figure \ref{sec:numres_param:subsec:parametric2:fig:normal_forces} depicts the evolution of normal forces over time for different parameter sets. We can clearly appreciate that the loss of prestress in the cables (a key feature in the study of leakage rates) strongly differs according to the pair of parameters studied. The observation of these quantities supports the choice of a logarithmic discretization for the construction of the parametric grid.

\begin{figure}[h!]
\begin{center}
\includegraphics[scale=0.85]{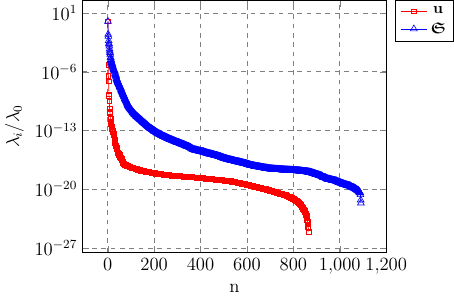} 
\end{center}
\caption{POD eigenvalues for the displacement and the generalized forces ($\bm{\mathfrak{S}}$) using a $\ell_2$ compression\label{sec:numres_param:subsec:parametric2:fig:eigenvales_ell2} for a parametric problem}
\end{figure}

Figure \ref{sec:numres_param:subsec:parametric2:fig:eigenvales_ell2} shows the decay of the POD eigenvalues when using the 25 HF snapshots. The decay is similar to that shown in Figure \ref{sec:numres:subsec:srpb:fig:eigenvales_ell2}. We notice that for the parametric case, the decay is fast and the gain in compression will be significant.

\begin{figure}[h!]
\begin{center}
\centering
\begin{subfigure}[b]{0.45\textwidth}
\begin{center}
\includegraphics[scale=0.85]{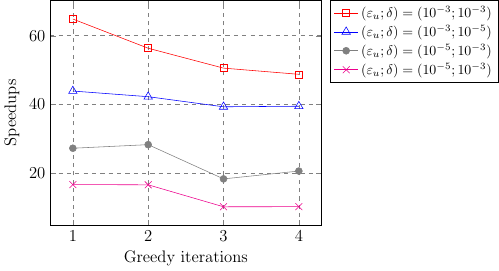} 
\end{center}
\caption{ \normalsize Speedups (POD on all HF snapshots)}
\end{subfigure}
\qquad
\centering
\begin{subfigure}[b]{0.45\textwidth}
\begin{center}
\includegraphics[scale=0.85]{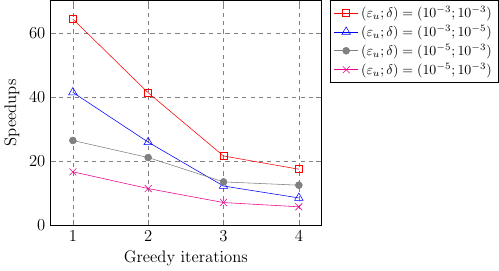} 
\end{center}
\caption{ \normalsize Speedups (Incremental POD)}
\end{subfigure}
\\
\centering
\begin{subfigure}[b]{0.45\textwidth}
\begin{center}
\includegraphics[scale=0.85]{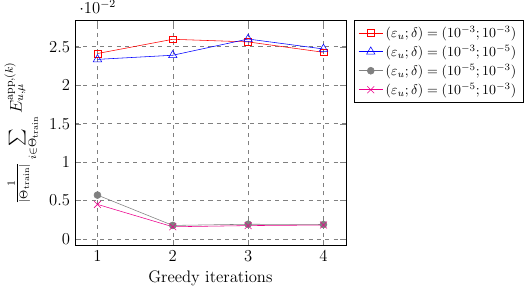} 
\end{center}
\caption{ \normalsize Average error (POD on all HF snapshots)}
\end{subfigure}
\qquad
\centering
\begin{subfigure}[b]{0.45\textwidth}
\begin{center}
\includegraphics[scale=0.85]{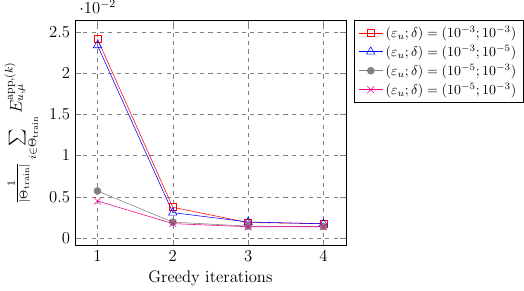} 
\end{center}
\caption{ \normalsize Average error (Incremental POD)}
\end{subfigure}
\end{center}
\caption{Speedups and average approximation errors on displacements fields for $\mu\in \Theta_{\rm train}$ using a training set of size $|\Theta_{\rm train}|=4$ for different compression tolerances ($\varepsilon$) and hyper-reduction parameters ($\delta$) and comparison between non-incremental and incremental POD \label{sec:numres_param:subsec:parametric2:speedups_avg_errors}}
\end{figure}

As a first test, we report a quick evaluation of the construction of a ROM on a smaller training set, consisting of 4 points. In other words, we take only the extremums of the 2d square to which all the parameters belong. The aim of this simpler case is to compare the two methodologies for building POD-reduced bases (in the parametric case) before presenting the case on the 25-point parametric case. Figure \ref{sec:numres_param:subsec:parametric2:speedups_avg_errors} depicts the speedups and approximation errors obtained after 4 iterations (the maximum number of iterations possible for this case) for different pairs of hyper-parameters used for ROM construction: number of modes and hyper-reduction tolerance. We observe that the hierarchical basis strategy leads to an increase in basis size (in our case), which reduces speedup and improves approximation quality (to below one percent). On the other hand, the use of full POD enables much better speedups to be maintained, while reducing the approximation error, but to a lesser extent. The same tradeoff applies to ROM construction as described above. In the case studied here, the regularity of the problem (at least for this set of parameters), prompts us to favor a POD on all snapshots (therefore, the basis is not hierarchical during iterations), in order to have the most efficient ROM both in terms of computational gain, while having reasonable approximation errors.

\begin{figure}[h!]
\begin{center}
\includegraphics[scale=0.85]{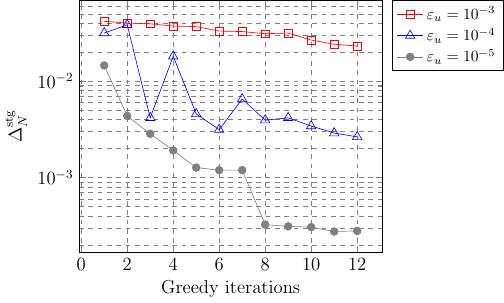} 
\end{center}
\caption{Maximum approximation error on unexplored parameters decreases during greedy iterations with an hyper-reduction parameter $\delta=10^{-5}$ \label{sec:numres_param:subsec:parametric2:fig:decreaseMaxError}}
\end{figure}

\begin{figure}[h!]
\begin{center}
\includegraphics[scale=0.65]{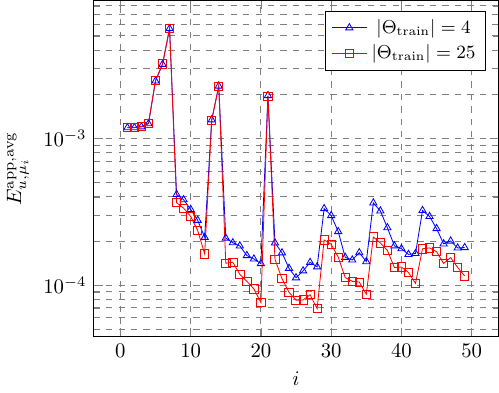} 
\end{center}
\caption{Statistical errors on the training set $\Theta_{\rm train}$, defined as a $5\times 5$ grid along the greedy iterations. Two strategies are compared: POD on all HF snapshots (\textcolor{red}{red}), and incremental POD (\textcolor{orange}{orange}) \label{sec:numres_param:subsec:parametric2:fig:statsErrorsTrain}}
\end{figure}

\begin{figure}[h!]
\begin{center}
\centering
\begin{subfigure}[b]{0.45\textwidth}
\includegraphics[scale=0.85]{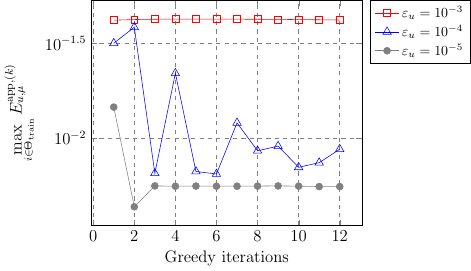} 
\caption{ \normalsize Maximum error\label{sec:numres_param:subsec:parametric2:fig:aemax}}
\end{subfigure}
\qquad
\centering
\begin{subfigure}[b]{0.45\textwidth}
\includegraphics[scale=0.85]{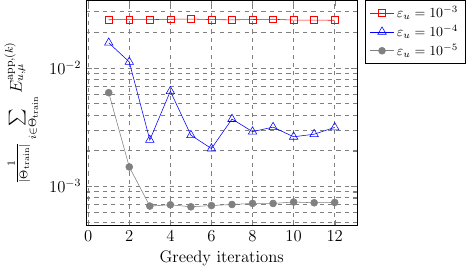} 
\caption{ \normalsize Average error\label{sec:numres_param:subsec:parametric2:fig:aeavg}}
\end{subfigure}
\end{center}
\caption{Average approximation errors on displacements fields for $\mu\in \Theta_{\rm train}$ using a training set of size $|\Theta_{\rm train}|=25$ and a non-incremental POD for different compression tolerances ($\varepsilon$) with an hyper-reduction parameter $\delta=10^{-5}$}
\end{figure}

\begin{figure}[h!]
\begin{center}
\begin{subfigure}[b]{0.45\textwidth}
\begin{center}
\includegraphics[scale=0.85]{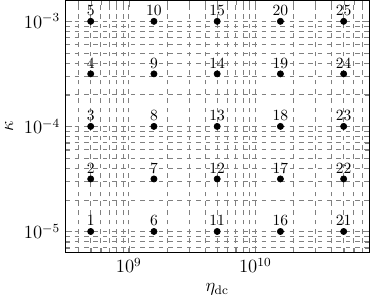} 
\end{center}
\caption{Indexes of parameters\label{sec:numres_param:subsec:parametric2:fig:avgae_allparams_greedy:sub:a}}
\end{subfigure}
\qquad
\centering
\begin{subfigure}[b]{0.45\textwidth}
\begin{center}
\includegraphics[scale=0.85]{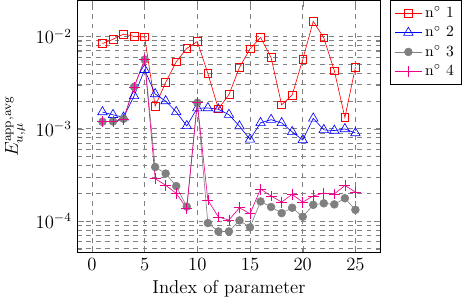} 
\end{center}
\caption{Approximation errors}
\end{subfigure}
\end{center}
\caption{Time-averaged approximation errors on displacement on the training set ($|\Theta_{\rm train}|=25$) for the first greedy iterations with an hyper-reduction parameter $\delta=10^{-5}$\label{sec:numres_param:subsec:parametric2:fig:avgae_allparams_greedy}} 
\end{figure}

Then, we apply this strategy to a larger training set ($|\Theta_{\rm train}| = 25$ parameters). Figure \ref{sec:numres_param:subsec:parametric2:fig:decreaseMaxError} represents the decay of the maximum approximation error on unexplored parameters (used to drive the greedy procedure). These successive choices clearly lead to a decrease in the maximum error (Figure \ref{sec:numres_param:subsec:parametric2:fig:aemax}) and the average error (Figure \ref{sec:numres_param:subsec:parametric2:fig:aeavg}) over the entire training set (explored and unexplored parameters). Scaling up for each parameter, Figure \ref{sec:numres_param:subsec:parametric2:fig:avgae_allparams_greedy} shows the time-averaged approximation errors for each parameter over the first iterations of the algorithm. As confirmed by the other figures, we observe that for the case studied, we have errors of the order of a few percent on all parameters (no more than ten percent) after just a few iterations. This is due to the relative regularity of the problem studied. Figure \ref{sec:numres_param:subsec:parametric2:fig:statsErrorsTrain} displays error statistics (median, quartiles) over the course of greedy iterations (5 by 5). We compare two approaches for incremental POD or POD on all snapshots, with error visualization, where we observe a decrease in medians over the iterations.

\subsubsection{Out-of-sample test for $\mathcal{P}\subset\R^2$}
All the above numerical results highlight the good approximation quality of the ROM on the training set. Nevertheless, it is crucial to further assess the methodology's suitability for out-of-sample parameters. To this end, we consider a 7-by-7 grid. This ensures that we get non-matching points. Then, we test the approximation quality of the ROM on this set, called the test set.

\begin{figure}[h!]
\begin{center}
\centering
\begin{subfigure}[b]{0.45\textwidth}
\begin{center}
\includegraphics[scale=0.65]{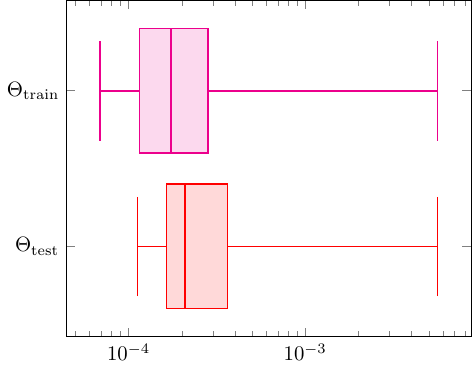} 
\end{center}
\caption{ \normalsize POD on all HF snapshots \label{sec:numres_param:subsec:parametric2:fig:TestVsTrain_5x5:sub:a}}
\end{subfigure}
\qquad
\centering
\begin{subfigure}[b]{0.45\textwidth}
\begin{center}
\includegraphics[scale=0.65]{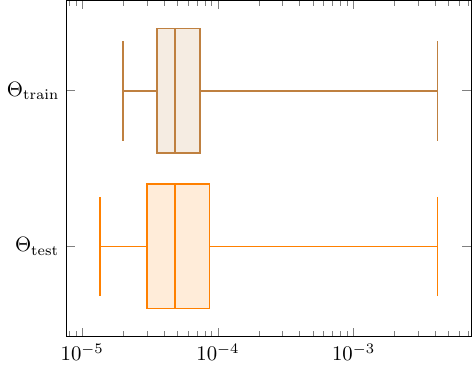} 
\end{center}
\caption{ \normalsize Incremental POD \label{sec:numres_param:subsec:parametric2:fig:TestVsTrain_5x5:sub:b}}
\end{subfigure}
\end{center}
\caption{Boxplot for a training set on a $5\times 5$ grid ($|\Theta_{\rm train}|=25$), verified on a test set on a $7\times 7$ grid ($|\Theta_{\rm test}|=49$). The quantities measured are the time-averaged errors on each set, for a ROM resulting from a greedy procedure, stopped after 5 iterations.\label{sec:numres_param:subsec:parametric2:fig:TestVsTrain_5x5}}
\end{figure}

Figure \ref{sec:numres_param:subsec:parametric2:fig:TestVsTrain_5x5} depicts boxplots for time-averaged approximation errors on the test set for the same training set for two sets of greedy strategies: one based on a POD on all snapshots (Figure \ref{sec:numres_param:subsec:parametric2:fig:TestVsTrain_5x5:sub:a}) and the other on an incremental POD (Figure \ref{sec:numres_param:subsec:parametric2:fig:TestVsTrain_5x5:sub:b}). From a statistical point of view, most of the test set features good approximation quality. The distribution of statistics across the two cases is consistent. For the POD on all snapshots, the error on the training set is of slightly higher quality than on the test set, while maintaining excellent approximation quality. Despite the simplicity of the case, it remains complex to perfectly capture the worst-case representations in the same way as the rest. Nevertheless, the worst-case error remains of the order of a few percent on the test set. For the case with incremental POD, the error quality between training and test sets is very similar, which is consistent with the fact that more modes are used than with POD on the snapshot set. Yet the difference between training and test sets is due to the smaller quartile spread on the training set (lower statistical dispersion), which is also coherent.

\begin{figure}[h!]
\begin{center}
\centering
\begin{subfigure}[b]{0.475\textwidth}
\begin{center}
\includegraphics[scale=0.65]{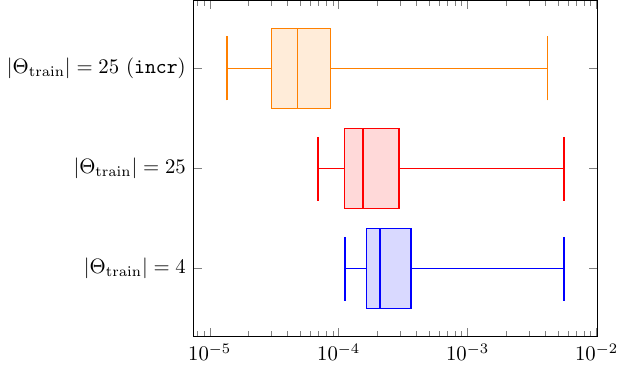} 
\end{center}
\caption{ \normalsize Boxplots \label{sec:numres_param:subsec:parametric2:fig:testErrorsAllStats:sub:a}}
\end{subfigure}
\qquad
\centering
\begin{subfigure}[b]{0.45\textwidth}
\begin{center}
\includegraphics[scale=0.6]{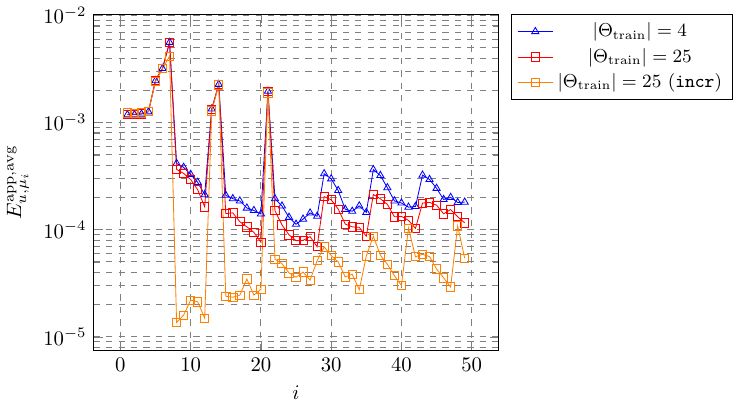} 
\end{center}
\caption{ \normalsize Approximation errors \label{sec:numres_param:subsec:parametric2:fig:testErrorsAllStats:sub:b}}
\end{subfigure}
\end{center}
\caption{Statiscal repartition of time-averaged errors generated by several ROMs on the same test set defined on a  $7\times 7$ grid ($|\Theta_{\rm test}|=49$). Three ROMs are compared (all obtained by a greedy process): built on a $2\times 2$ training grid with POD on all HF snapshots (\textcolor{blue}{blue}), on a $5\times 5$ training grid with POD on all HF snapshots (\textcolor{red}{red}), and on a $5\times 5$ training grid with an incremental POD (\textcolor{orange}{orange}).  Figure \ref{sec:numres_param:subsec:parametric2:fig:testErrorsAllStats:sub:a} is a boxplot of time-averaged errors on $\Theta_{\rm test}$ and \ref{sec:numres_param:subsec:parametric2:fig:testErrorsAllStats:sub:b} is the time-averaged errors according to the number of the parameters in the $\Theta_{\rm test}$ (numerotation is similar to Figure \ref{sec:numres_param:subsec:parametric2:fig:avgae_allparams_greedy:sub:a}, but on a $7\times 7$ grid)\label{sec:numres_param:subsec:parametric2:fig:testErrorsAllStats}}
\end{figure}

In a second step, we can also compare the greedy approaches with each other in terms of their behavior on the test set (Figure \ref{sec:numres_param:subsec:parametric2:fig:testErrorsAllStats}). As can be expected, the poorest approximation case matches the case with the smallest training set size, followed by the case with 25 points and total POD, followed by a case with 25 points and incremental POD. This analysis is reflected in the boxplots (see Figure \ref{sec:numres_param:subsec:parametric2:fig:testErrorsAllStats:sub:a}), as well as in the plot of errors as a function of parameter indices (indices are distributed in a similar way to discretization on a 5x5 grid). 

\section{Conclusion}

We proposed and validated a methodology for the construction of ROMs for multi-modeling problems, with an application to a standard section of prestressed concrete NCB. This involves several aspects. First, we devised a robust numerical method, suitable for use with industrially-constrained codes, providing ROMs designed to replicate the behavior of prestressed concrete with high speedups and good approximation errors. Furthermore, we proposed an adaptive approach to iteratively enrich the reduced model on a set of parameters. These two points are presented theoretically and validated numerically. Second, we have also succeeded in producing a ROM that can be used for real engineering applications, in that it provides a good representation of the variables of interest used in practice by engineers, whether for structural state analysis (leakage rate study) or for in-depth data analysis (data-assimilation problem, Bayesian approaches). \\

Much work is currently underway to make further progress in several directions. First, these  promising results are validated on fairly coarse meshes (although used in practice) and on smaller parametric spaces. Efforts are currently underway to evaluate these approaches by increasing the dimension of the parametric vector, and of the snapshot vectors considered (mesh refinement). Second, the approach adopted is a strong greedy process and relies on comparison with known HF snapshots. This leads to significant offline computation costs, since it requires \emph{a priori} knowledge of these solutions. This is a particular limitation when scaling up. Previous efforts have focused on the construction of low-cost a posteriori error indicators. The efficiency of these indicators in steering greedy search (within a weak-greedy context) has been demonstrated for problems featuring internal variables. The problems presented here are somewhat more intricate from a theoretical standpoint in mechanics (THM and multimodeling), and consequently, pose challenges for robust implementation in an industrial-grade FE code. Ongoing efforts are being made  to broaden the application of indicators to tackle these challenges. Research is also underway to make this methodology still applicable when the parameter space becomes larger. Finally, the coupling of the ROM methodology with the optimization problems mentioned above, and in particular data assimilation, in order to reduce the resolution time, are being studied.

\section*{Acknowledgements}

This work was supported by ANRT (French National Association for Research and Technology) and EDF. We extend our gratitude to the \texttt{code$\_$aster }development team and all contributors to the code. Our focus has been on utilizing and advancing the Python library Mordicus, supported by a 'French Fonds Unique Interministériel' (FUI) project and designed as a tool for the advancement of model reduction methods tailored for industrial applications.

\appendix

\end{document}